# The linear Cahn-Hilliard equation with an interface

Andreas Chatziafratis [1,*], Alain Miranville [2], Tohru Ozawa [3]

[1] Department of Mathematics, National and Kapodistrian University of Athens, Greece

[2] Henan Normal University, School of Mathematics and Statistics, Xinxiang, China
Laboratoire de Mathématiques Appliquées du Havre (LMAH), Université Le Havre Normandie, France

[3] Department of Applied Physics, Waseda University, Tokyo, Japan

May 18, 2026

**Abstract**. We obtain new integral representations, expressed as contour integrals in the complex Fourier plane, for the solution of fully nonhomogeneous interface problems for the linearized Cahn–Hilliard equation with arbitrary initial data on the line and general interface conditions prescribed at the origin. Cahn–Hilliard-type models emerge in applied mathematics in connection to a spectacular variety of phenomena of mathematical physics, continuum mechanics, chemistry and biology. A novel implementation of Fokas' unified transform method is in force herein for a fourth-order operator for the first time, with particular challenges arising due to the nature and the generality of the problems under consideration. Our explicit formulae directly lend themselves to exploration of the solution's qualitative properties such as regularity and asymptotic behavior. This work is also useful in the investigation of well-posedness for nonlinear counterparts as well as in the study of free-boundary and diffuse-interface problems.

## 1. Introduction

The Cahn–Hilliard equation [1,2] constitutes a canonical model in the study of phase separation dynamics and diffuse interface evolution in binary mixtures. Since its introduction, it has served as a fundamental example of a mass-conserving fourth-order parabolic equation, exhibiting a rich interplay between dissipation, interfacial structure, and nonlinear pattern formation. Owing to these features, the equation has generated extensive interest in both the analysis of nonlinear partial differential equations and the mathematical theory of complex materials; see [3-6] and references therein.

Interface formulations of Cahn–Hilliard-type equations arise naturally in contemporary continuum modelling. In materials science, transmission conditions across interfaces are used to model heterogeneous media, composite alloys, grain boundaries, and thin multilayer structures, where distinct material phases interact through localized coupling conditions [3-9]. Closely related models appear in the theory of lithium-ion batteries, e.g. [10,11], particularly in phase-separating electrode materials such as lithium iron phosphate, where Cahn–Hilliard dynamics govern ion concentration evolution while interfaces represent material discontinuities or internal boundaries between active components. Moreover, diffuse-interface models based on the Cahn–Hilliard framework have become increasingly important in biological fluid mechanics and soft matter theory, e.g. [12-17], including the modelling of lipid membranes, tumour growth, intracellular phase separation, and multiphase biological flows. In such settings, localised interface conditions may encode membrane permeability, transport constraints, or coupling between distinct physical compartments. Notably, in addition, the linear Cahn–Hilliard operator may be viewed as a perturbation of the biharmonic heat operator, the latter arising in surface diffusion phenomena, thin-film dynamics, and thermal grooving models [18]. These developments further motivate the analysis of explicit interface problems for higher-order phase-field equations on unbounded domains.

In the present work, we investigate interface problems for the one-dimensional linearized Cahn–Hilliard equation posed on the real line supplemented with interface conditions prescribed at the origin ( $x = 0$ ). The interface conditions couple the solutions on the two half-lines through constraints involving traces of the solution and its derivatives, thereby modelling transmission phenomena across a localized defect or material discontinuity. Such problems arise naturally in heterogeneous media, phase-field models with localized interactions, and in the analysis of composite systems where discontinuous constitutive behavior is present at an interface.

---

*corresponding author, e-mail: chatziafrati@math.uoa.gr

From an analytical perspective, interface problems for higher-order evolutionary equations present substantial difficulties. In contrast to second-order diffusive equations, the fourth-order structure of the Cahn–Hilliard operator requires the prescription and compatibility of multiple boundary or transmission conditions, while simultaneously preserving the underlying conservation and dissipative properties of the flow. Moreover, on unbounded domains the interaction between the interface and the dispersive-diffusive character of the linearized dynamics leads to a nontrivial spectral structure, making explicit solution representations considerably more delicate.

The principal aim of this paper is to derive explicit solution formulae for a broad class of one-dimensional (representing a symmetry in higher-dimensional configurations) interface problems associated with the linear Cahn–Hilliard equation. Our approach combines the Fokas' unified transform methodology [19-23] with complex-analytic techniques in order to obtain integral representations of the solution directly in terms of the prescribed initial, interface and forcing data. The Fokas method has so far been successfully implemented to interface problems for second- and third-order linear evolution PDE; see e.g. [24,25]. In recent works, the UTM was effectively extended for stationary contact problems [26] and for mixed-derivative nonlocal dispersive PDE [27].

Here, for the first time, an analogous study for a *fourth-order* PDE is undertaken. In particular, the analysis yields exact characterisations of the solution operators and clarifies the role played by the interface conditions in determining the qualitative behavior of the associated evolution. The presence of a fourth-order operator, together with the lower-order term and general form of interface conditions that we address, introduces structural features and analytical complications absent from the existing literature.

Beyond their intrinsic interest, the explicit formulae obtained in this work provide a useful framework for further investigations of nonlinear Cahn–Hilliard interface problems. In particular, they furnish a natural starting point for studying asymptotic behavior, well-posedness, stability questions and so forth (as has been done, for instance, in [5,28-39]), as well as for analyzing even more general (e.g. higher-order, nonlinear, boundary-layer) evolutionary systems posed in heterogeneous geometries, problems to be considered elsewhere.

Let us now therefore consider the following:

***Problem*** Solve

$$(1.1)\qquad \begin{cases} u_t^R = \alpha_R u_{xx}^R - \beta_R u_{xxxx}^R + f_R \\ u^R(x,0) = u_0^R(x) \end{cases}$$

for $u^R(x,t)$, with $x>0$, $t>0$, and

$$(1.2)\qquad \begin{cases} u_t^L = \alpha_R u_{xx}^L - \beta_R u_{xxxx}^L + f_L \\ u^L(x,0) = u_0^L(x) \end{cases}$$

for $u^L(x,t)$, with $x<0$, $t>0$, such that

$$(1.3)\begin{cases} \gamma_{11}u^R(0,t)+\gamma_{12}u_x^R(0,t)+\gamma_{13}u_{xx}^R(0,t)+\gamma_{14}u_{xxx}^R(0,t)+\gamma_{15}u^L(0,t)+\gamma_{16}u_x^L(0,t)+\gamma_{17}u_{xx}^L(0,t)+\gamma_{18}u_{xxx}^L(0,t)=I_1(t) \\ \gamma_{21}u^R(0,t)+\gamma_{22}u_x^R(0,t)+\gamma_{23}u_{xx}^R(0,t)+\gamma_{24}u_{xxx}^R(0,t)+\gamma_{25}u^L(0,t)+\gamma_{26}u_x^L(0,t)+\gamma_{27}u_{xx}^L(0,t)+\gamma_{28}u_{xxx}^L(0,t)=I_2(t) \\ \gamma_{31}u^R(0,t)+\gamma_{32}u_x^R(0,t)+\gamma_{33}u_{xx}^R(0,t)+\gamma_{34}u_{xxx}^R(0,t)+\gamma_{35}u^L(0,t)+\gamma_{36}u_x^L(0,t)+\gamma_{37}u_{xx}^L(0,t)+\gamma_{38}u_{xxx}^L(0,t)=I_3(t) \\ \gamma_{41}u^R(0,t)+\gamma_{42}u_x^R(0,t)+\gamma_{43}u_{xx}^R(0,t)+\gamma_{44}u_{xxx}^R(0,t)+\gamma_{45}u^L(0,t)+\gamma_{46}u_x^L(0,t)+\gamma_{47}u_{xx}^L(0,t)+\gamma_{48}u_{xxx}^L(0,t)=I_4(t) \end{cases}$$

for $t>0$.

***Data*** We are given the initial values $u_0^R(x)$ ($x\ge 0$) and $u_0^L(x)$ ($x\le 0$), and the forcing functions $f_R(t,x)$ ($x\ge 0$, $t\ge 0$) and $f_L(t,x)$, ($x\le 0$, $t\ge 0$), satisfying the following assumptions:

(1.4) $u_0^R\in\mathcal{S}([0,\infty))$, $u_0^L\in\mathcal{S}((-\infty,0])$,

(1.5) $f_R(t,\cdot)\in\mathcal{S}([0,\infty))$, $f_L(t,\cdot)\in\mathcal{S}((-\infty,0])$, uniformly for $t$ in compact subsets of $[0,\infty)$.

Assumptions on the parameters:

(1.6) $\alpha_R>0$, $\beta_R>0$, $\alpha_L>0$, $\beta_L>0$, $\gamma_{jk}\in\mathbb{R}$.

In our Theorem, stated in *Section 4*, more specific assumptions on the constants $\gamma_{jk}$ are made.

## 2. Derivation of the integral representation for $x > 0$

Setting $\omega_R(\lambda) = \alpha_R\lambda^2 + \beta_R\lambda^4$, we write the equation $u_t^R = \alpha_R u_{xx}^R - \beta_R u_{xxxx}^R + f_R$ in the form:

$$\frac{\partial}{\partial t}[e^{-i\lambda x+\omega_R(\lambda)t}u^R(x,t)] - \frac{\partial}{\partial x}\{e^{-i\lambda x+\omega_R(\lambda)t}[-\beta_R u_{xxx}^R - \beta_R i\lambda u_{xx}^R + (\beta_R\lambda^2+\alpha_R)u_x^R + i\lambda(\beta_R\lambda^2+\alpha_R)u^R]\}$$
$$= e^{-i\lambda x+\omega_R(\lambda)t} f_R(x,t),$$

and Green's theorem gives the global relation:

$$(2.1)\quad \hat{u}^R(\lambda,t)e^{\omega_R(\lambda)t} = \hat{u}_0^R(\lambda) - i\lambda(\beta_R\lambda^2+\alpha_R)\tilde{g}_{0,R}(\omega_R(\lambda),t) - (\beta_R\lambda^2+\alpha_R)\tilde{g}_{1,R}(\omega_R(\lambda),t)$$
$$+\beta_R i\lambda\tilde{g}_{2,R}(\omega_R(\lambda),t) + \beta_R\tilde{g}_{3,R}(\omega_R(\lambda),t) + \tilde{\hat{f}}_R(\lambda,\omega_R(\lambda),t),\ \operatorname{Im}\lambda \le 0,$$

where we have set $g_{0,R}(t) = u^R(0,t)$, $g_{1,R}(t) = u_x^R(0,t)$, $g_{2,R}(t) = u_{xx}^R(0,t)$, $g_{3,R}(t) = u_{xxx}^R(0,t)$ and

$$\tilde{g}_{j,R}(\omega_R(\lambda),t) = \int_{\tau=0}^{t} g_{j,R}(\tau)e^{i\omega_R(\lambda)\tau}d\tau,\ j = 0,1,2,3.$$

Multiplying (2.1) by $e^{i\lambda x-\omega_R(\lambda)t}$ and integrating, we obtain

$$(2.2)\quad 2\pi u^R(x,t) = \int_{-\infty}^{\infty} e^{i\lambda x-\omega_R(\lambda)t}\hat{u}_0^R(\lambda)d\lambda - i\int_{-\infty}^{\infty}\lambda e^{i\lambda x-\omega_R(\lambda)t}(\beta_R\lambda^2+\alpha_R)\tilde{g}_0(\omega_R(\lambda),t)d\lambda$$
$$-\int_{-\infty}^{\infty} e^{i\lambda x-\omega_R(\lambda)t}(\beta_R\lambda^2+\alpha_R)\tilde{g}_{1,R}(\omega_R(\lambda),t)d\lambda + \beta_R i\int_{-\infty}^{\infty} e^{i\lambda x-\omega_R(\lambda)t}\lambda\tilde{g}_{2,R}(\omega_R(\lambda),t)d\lambda$$
$$+\beta_R\int_{-\infty}^{\infty} e^{i\lambda x-\omega_R(\lambda)t}\tilde{g}_{3,R}(\omega_R(\lambda),t)d\lambda + \int_{-\infty}^{\infty} e^{i\lambda x-\omega_R(\lambda)t}\tilde{\hat{f}}_R(\lambda,\omega_R(\lambda),t)d\lambda,\ x>0,\ t>0.$$

Let us consider the set

$$\Omega_R = \{\lambda : \operatorname{Im}\lambda \ge 0\ \&\ \operatorname{Re}\omega_R(\lambda) = \operatorname{Re}(\alpha_R\lambda^2 + \beta_R\lambda^4) \le 0\}.\ \text{(See fig1.)}$$

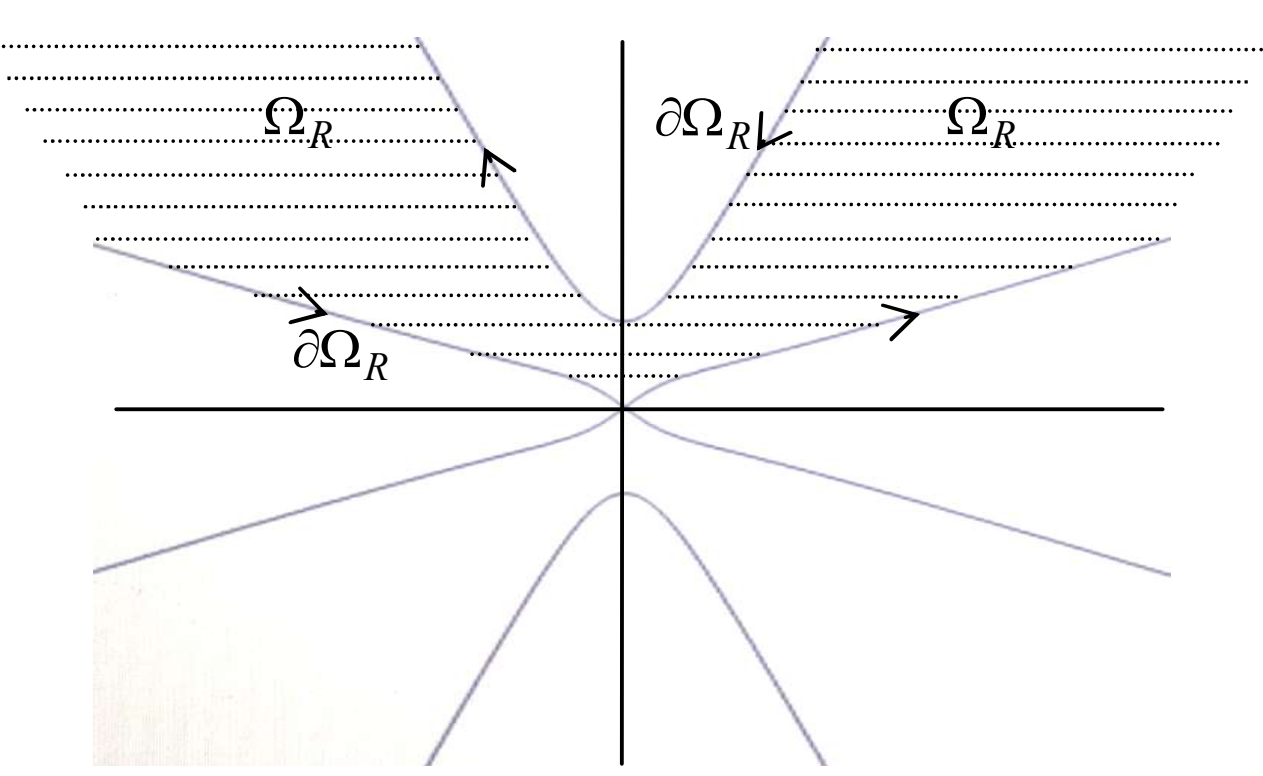

**Fig1** *The set* $\Omega_R$ *and its boundary* $\partial\Omega_R$.

Then we write (2.2) as follows:

$$(2.3)\quad 2\pi u^R(x,t) = \int_{-\infty}^{\infty} e^{i\lambda x-\omega_R(\lambda)t}\hat{u}_0^R(\lambda)d\lambda - i\int_{\partial\Omega_R}\lambda e^{i\lambda x-\omega_R(\lambda)t}(\beta_R\lambda^2+\alpha_R)\tilde{g}_0(\omega_R(\lambda),t)d\lambda$$

$$-\int_{\partial\Omega_R} e^{i\lambda x-\omega_R(\lambda)t}(\beta_R\lambda^2+\alpha_R)\widetilde{g}_{1,R}(\omega_R(\lambda),t)d\lambda+\beta_R i\int_{\partial\Omega_R} e^{i\lambda x-\omega_R(\lambda)t}\lambda\widetilde{g}_{2,R}(\omega_R(\lambda),t)d\lambda$$

$$+\beta_R\int_{\partial\Omega_R} e^{i\lambda x-\omega_R(\lambda)t}\widetilde{g}_{3,R}(\omega_R(\lambda),t)d\lambda+\int_{-\infty}^{\infty} e^{i\lambda x-\omega_R(\lambda)t}\widetilde{\hat{f}}_R(\lambda,\omega_R(\lambda),t)d\lambda\,,\ x>0\,,\ t>0\,.$$

Next we consider the mappings $\mu=\psi_R(\lambda)$ and $\lambda=\varphi_R(\mu)$ so that

$$\omega_R(\lambda)=\alpha_R\lambda^2+\beta_R\lambda^4=\mu^4=:\omega_0(\mu)\,,$$

which are defined as follows:

$$\mu=\psi_R(\lambda)=\begin{cases}\sqrt[4]{\alpha_R\lambda^2+\beta_R\lambda^4}\ \text{if}\ \operatorname{Re}\lambda>0,\\ i\sqrt[4]{\alpha_R\lambda^2+\beta_R\lambda^4}\ \text{if}\ \operatorname{Re}\lambda<0\,.\end{cases}\quad\text{and}\ \lambda=\varphi_R(\mu)=\begin{cases}\sqrt{\dfrac{-\alpha_R+\sqrt{4\beta_R\mu^4+\alpha_R^2}}{2\beta_R}}\ \text{if}\ \operatorname{Re}\mu>0,\\ \sqrt{\dfrac{-\alpha_R-\sqrt{4\beta_R\mu^4+\alpha_R^2}}{2\beta_R}}\ \text{if}\ \operatorname{Re}\mu<0\,.\end{cases}$$

(The choices of the functions $\sqrt[4]{\cdots}$ and $\sqrt{\cdots}$ are described in fig 2 and fig 3.)

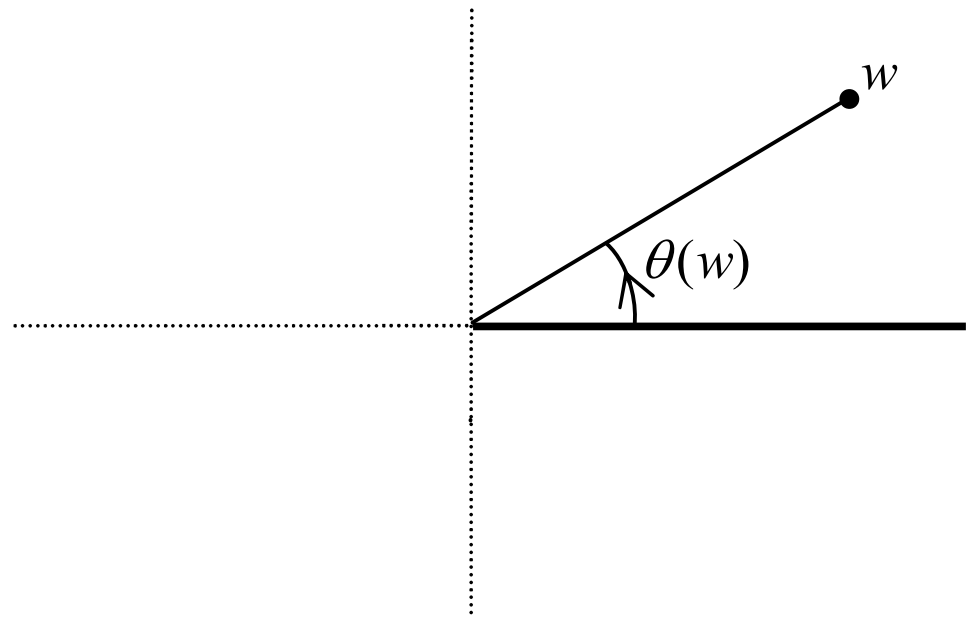

**Fig.2** *The choice of the angle $\theta(w)$ in the definition of the 4th root function:*
$\sqrt[4]{w}=\sqrt[4]{|w|}e^{i\theta(w)/4}:\ 0\le\theta(w)<2\pi$ .

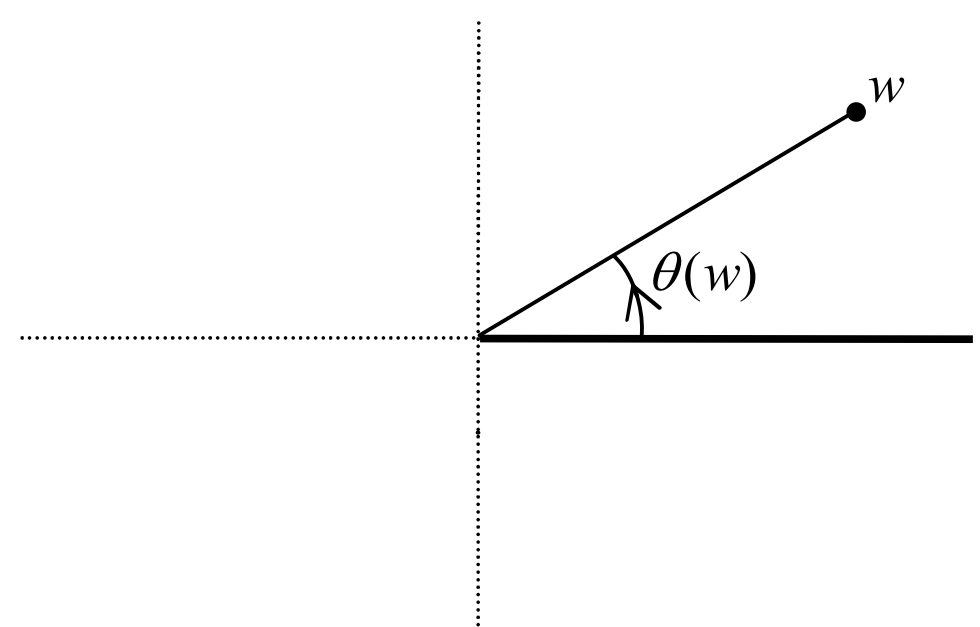

**Fig.3** *The choice of the angle $\theta(w)$ in the definition of the square root function:*
$\sqrt{w}=\sqrt{|w|}e^{i\theta(w)/2}:\ 0\le\theta(w)<2\pi$ .

Then $\psi_R=\psi_R(\lambda)$ is a biholomorphic mapping from an open neighborhood of $\Omega_R-[0,i\sqrt{\alpha_R/\beta_R}]$ to an open neighborhood of the set $\Delta_R:=\{\mu\in\mathbb{C}-\{0\}:\frac{\pi}{8}\le\arg\mu\le\frac{3\pi}{8}\ or\ \frac{5\pi}{8}\le\arg\mu\le\frac{7\pi}{8}\}$, with inverse $\varphi_R=\varphi_R(\mu)$.

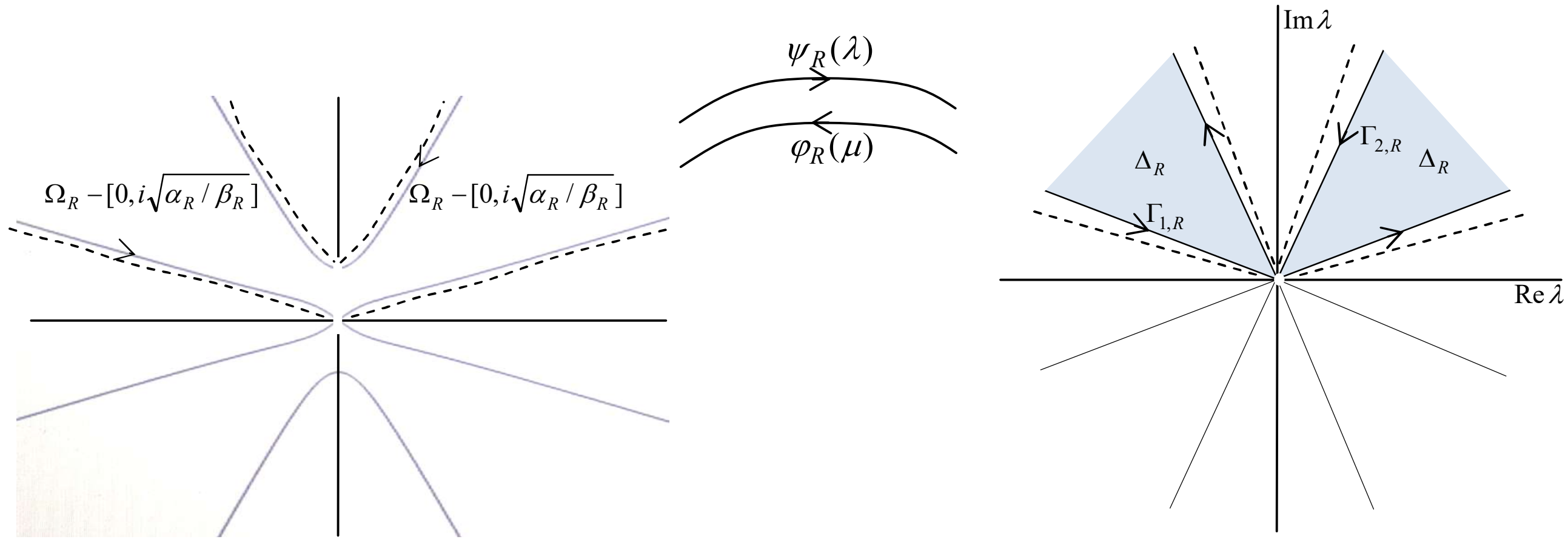

**Fig.4** *The mappings* $\mu = \psi_R(\lambda)$ *and* $\lambda = \varphi_R(\mu)$.

Setting $\lambda = \varphi_R(\mu)$ and changing the variable from $\lambda$ to $\mu$, in the integrals over $\partial\Omega_R$, we write (2.3) in the following way:

$$(2.4)\quad 2\pi u^R(x,t) = \int_{-\infty}^{\infty} e^{i\lambda x - \omega_R(\lambda)t} \hat{u}_0^R(\lambda) d\lambda - i \int_{\Gamma_R} \varphi_R(\mu) e^{i\varphi_R(\mu)x - \omega_0(\mu)t} (\beta_R \varphi_R^2(\mu) + \alpha_R) \tilde{g}_{0,R}(\omega_0(\mu),t) \frac{d\varphi_R(\mu)}{d\mu} d\mu$$

$$- \int_{\Gamma_R} e^{i\varphi_R(\mu)x - \omega_0(\mu)t} (\beta_R \varphi_R^2(\mu) + \alpha_R) \tilde{g}_{1,R}(\omega_0(\mu),t) \frac{d\varphi_R(\mu)}{d\mu} d\mu$$

$$+ \beta_R i \int_{\Gamma_R} \varphi_R(\mu) e^{i\varphi_R(\mu)x - \omega_0(\mu)t} \tilde{g}_{2,R}(\omega_0(\mu),t) \frac{d\varphi_R(\mu)}{d\mu} d\mu$$

$$+ \beta_R \int_{\Gamma_R} e^{i\varphi_R(\mu)x - \omega_0(\mu)t} \tilde{g}_{3,R}(\omega_0(\mu),t) \frac{d\varphi_R(\mu)}{d\mu} d\mu + \int_{-\infty}^{\infty} e^{i\lambda x - \omega_R(\lambda)t} \tilde{\hat{f}}_R(\lambda, \omega_R(\lambda), t) d\lambda,\ x > 0,\ t > 0,$$

where $\Gamma_R := \partial\Delta_R = \Gamma_{1,R} + \Gamma_{2,R}$ with $\Gamma_{1,R} := \{\lambda : \arg\lambda = \frac{\pi}{8}\ or\ \frac{7\pi}{8}\}$ and $\Gamma_{2,R} := \{\lambda : \arg\lambda = \frac{3\pi}{8}\ or\ \frac{7\pi}{8}\}$. (See fig 5.)

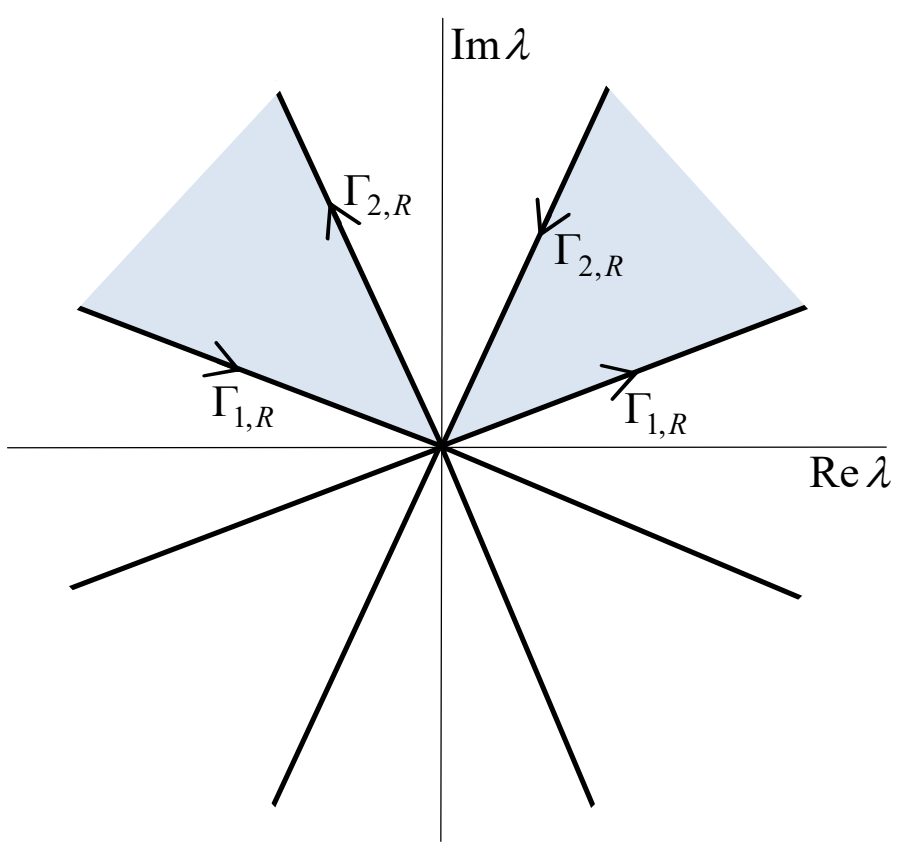

**Fig.5** *The contours* $\Gamma_{1,R}$ *and* $\Gamma_{2,R}$, *and their orientation.*

***Remarks*** **(1)** We point out that

$$\lim_{\substack{\mu\to 0 \\ \arg\mu = \pi/8\ or\ \arg\mu = 7\pi/8}} \varphi_R(\mu) = 0 \quad \text{and} \quad \lim_{\substack{\mu\to 0 \\ \arg\mu = 3\pi/8\ or\ \arg\mu = 5\pi/8}} \varphi_R(\mu) = i\sqrt{\alpha_R/\beta_R}.$$

**(2)** The derivative

$$\frac{d\lambda}{d\mu}=\frac{d\varphi_R(\mu)}{d\mu}=\begin{cases}16\beta_R\mu^3\left(\dfrac{-\alpha_R+\sqrt{4\beta_R\mu^4+\alpha_R^2}}{2\beta_R}\right)^{-1/2}(4\beta_R\mu^4+\alpha_R^2)^{-1/2} & \text{if } \operatorname{Re}\mu>0\\ 16\beta_R i\mu^3\left(\dfrac{-\alpha_R+\sqrt{4\beta_R\mu^4+\alpha_R^2}}{2\beta_R}\right)^{-1/2}(4\beta_R\mu^4+\alpha_R^2)^{-1/2} & \text{if } \operatorname{Re}\mu<0,\end{cases}$$

whence

(2.5) $\dfrac{d\varphi_R(\mu)}{d\mu}=\mathrm{O}(\mu^2)$, as $\mu\to 0$ with $\mu\in\Delta_R$ and $\dfrac{d\varphi_R(\mu)}{d\mu}=\mathrm{O}(1)$, as $\mu\to\infty$ with $\mu\in\Delta_R$.

In order to compute the integrals which contain the unknown quantities $\widetilde{g}_{j,R}(\omega_0(\mu),t)$, $j=0,1,2,3$, we need to observe that if $\rho_R(\lambda)=\sqrt{\lambda^2+\dfrac{\alpha_R}{\beta_R}}$ then

$$\omega_R(-\lambda)=\omega_R(i\rho_R(\lambda))=\omega_R(-i\rho_R(\lambda))=\omega_R(\lambda)\,,\ \forall\lambda\,.$$

We choose the following holomorphic branch of the square root of $\lambda^2+\dfrac{\alpha_R}{\beta_R}$, $\rho_R:\Theta_R\to\mathbb{C}$, with

$$\rho_R(\lambda):=\sqrt{\left|\lambda^2+\frac{\alpha_R}{\beta_R}\right|}\exp\{i[\theta_1(\lambda)+\theta_2(\lambda)+\pi]/2\}\ \text{(see fig 6)}$$

and

$$\Theta_R:=\mathbb{C}-\{(-i\infty,-\sqrt{\tfrac{\alpha_R}{\beta_R}}i]\cup[\sqrt{\tfrac{\alpha_R}{\beta_R}}i,i\infty)\}\,.$$

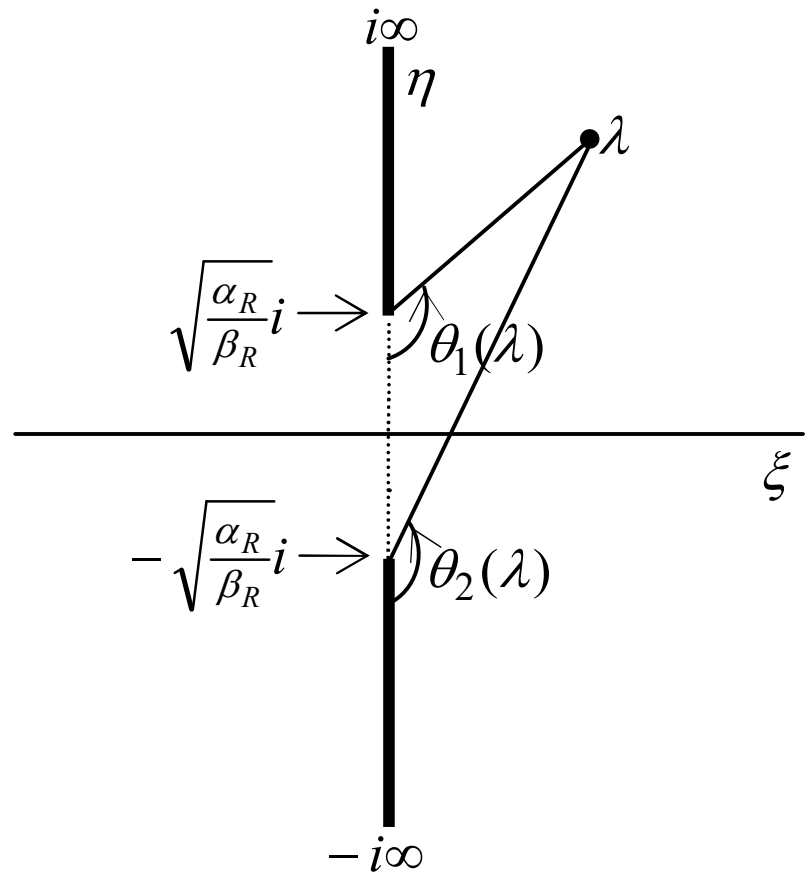

**Fig. 6** With $-\pi<\theta_1(\lambda)<\pi$ and $0<\theta_2(\lambda)<2\pi$, the function $\rho_R(\lambda)$ is analytic for $\lambda\in\Theta_R$.

Thus, the global relation (2.1) gives

$$\hat{u}^R(-\lambda,t)e^{\omega_R(\lambda)t}=\hat{u}_0^R(-\lambda)+i\lambda(\beta_R\lambda^2+\alpha_R)\widetilde{g}_{0,R}(\omega_R(\lambda),t)-(\beta_R\lambda^2+\alpha_R)\widetilde{g}_{1,R}(\omega_R(\lambda),t)$$
$$-\beta_R i\lambda\widetilde{g}_{2,R}(\omega_R(\lambda),t)+\beta_R\widetilde{g}_{3,R}(\omega_R(\lambda),t)+\widetilde{\hat{f}}(-\lambda,\omega_R(\lambda),t)\,,\ \operatorname{Im}\lambda\ge 0\,,$$

and

$$\hat{u}_R(i\rho_R(\lambda),t)e^{\omega_R(\lambda)t}=\hat{u}_0^R(i\rho_R(\lambda))-\beta_R\lambda^2\rho_R(\lambda)\widetilde{g}_{0,R}(\omega_R(\lambda),t)+\beta_R\lambda^2\widetilde{g}_{1,R}(\omega_R(\lambda),t)$$
$$-\beta_R\rho_R(\lambda)\widetilde{g}_{2,R}(\omega_R(\lambda),t)+\beta_R\widetilde{g}_{3,R}(\omega_R(\lambda),t)+\widetilde{\hat{f}}_R(i\rho_R(\lambda),\omega_R(\lambda),t)\,,\ \operatorname{Im}[i\rho_R(\lambda)]\le 0\,.$$

Setting $\lambda=\varphi_R(\mu)$ we obtain

(2.6) $\hat{u}^R(-\varphi_R(\mu),t)e^{\omega_0(\mu)t} = \hat{u}_0^R(-\varphi_R(\mu)) + i\varphi_R(\mu)(\beta_R\varphi_R^2(\mu)+\alpha_R)\widetilde{g}_{0,R}(\omega_0(\mu),t) - (\beta_R\varphi_R^2(\mu)+\alpha_R)\widetilde{g}_{1,R}(\omega_0(\mu),t)$

$$- \beta_R i\varphi_R(\mu)\widetilde{g}_{2,R}(\omega_0(\mu),t) + \beta_R\widetilde{g}_{3,R}(\omega_0(\mu),t) + \widetilde{\widetilde{f}}(-\varphi_R(\mu),\omega_0(\mu),t)\,,\ \mu\in\Delta_R\,,$$

and

(2.7) $\hat{u}_R(i\sigma_R(\mu),t)e^{\omega_0(\mu)t} = \hat{u}_0^R(i\sigma_R(\mu)) - \beta_R\varphi_R^2(\mu)\sigma_R(\mu)\widetilde{g}_{0,R}(\omega_0(\mu),t) + \beta_R\varphi_R^2(\mu)\widetilde{g}_{1,R}(\omega_0(\mu),t)$

$$- \beta_R\sigma_R(\mu)\widetilde{g}_{2,R}(\omega_0(\mu),t) + \beta_R\widetilde{g}_{3,R}(\omega_0(\mu),t) + \widetilde{\widetilde{f}}_R(i\sigma_R(\mu),\omega_0(\mu),t)\,,\ \mu\in\Delta_R\,,$$

where $\sigma_R(\mu) := \rho_R(\varphi_R(\mu))$.

## 3. Derivation of the integral representation for $x<0$

Working as in Section 2, we set $\omega_L(\lambda) = \alpha_L\lambda^2 + \beta_L\lambda^4$, we write the equation $u_t^L = \alpha_L u_{xx}^L - \beta_L u_{xxxx}^L + f_L$ in the form:

$$\frac{\partial}{\partial t}[e^{-i\lambda x+\omega_L(\lambda)t}u^L(x,t)] - \frac{\partial}{\partial x}\{e^{-i\lambda x+\omega_L(\lambda)t}[-\beta_L u_{xxx}^L - \beta_L i\lambda u_{xx}^L + (\beta_L\lambda^2+\alpha_L)u_x^L + i\lambda(\beta_L\lambda^2+\alpha_L)u^L]\}$$
$$= e^{-i\lambda x+\omega_L(\lambda)t}f_L(x,t)\,,$$

and Green's theorem gives the global relation:

(3.1) $\hat{u}^L(\lambda,t)e^{\omega_L(\lambda)t} = \hat{u}_0^L(\lambda) + i\lambda(\beta_L\lambda^2+\alpha_L)\widetilde{g}_{0,L}(\omega_L(\lambda),t) + (\beta_L\lambda^2+\alpha_L)\widetilde{g}_{1,L}(\omega_L(\lambda),t)$

$$- \beta_L i\lambda\widetilde{g}_{2,L}(\omega_L(\lambda),t) - \beta_L\widetilde{g}_{3,L}(\omega_L(\lambda),t) + \widetilde{\widetilde{f}}_L(\lambda,\omega_L(\lambda),t)\,,\ \operatorname{Im}\lambda\geq 0\,.$$

Therefore,

(3.2) $2\pi u^L(x,t) = \int_{-\infty}^{\infty} e^{i\lambda x-\omega_L(\lambda)t}\hat{u}_0^L(\lambda)d\lambda + i\int_{-\infty}^{\infty}\lambda e^{i\lambda x-\omega_L(\lambda)t}(\beta_L\lambda^2+\alpha_L)\widetilde{g}_0(\omega_L(\lambda),t)d\lambda$

$$+ \int_{-\infty}^{\infty} e^{i\lambda x-\omega_L(\lambda)t}(\beta_L\lambda^2+\alpha_L)\widetilde{g}_{1,L}(\omega_L(\lambda),t)d\lambda - \beta_L i\int_{-\infty}^{\infty} e^{i\lambda x-\omega_L(\lambda)t}\lambda\widetilde{g}_{2,L}(\omega_L(\lambda),t)d\lambda$$

$$- \beta_L\int_{-\infty}^{\infty} e^{i\lambda x-\omega_L(\lambda)t}\widetilde{g}_{3,L}(\omega_L(\lambda),t)d\lambda + \int_{-\infty}^{\infty} e^{i\lambda x-\omega_L(\lambda)t}\widetilde{\widetilde{f}}_L(\lambda,\omega_L(\lambda),t)d\lambda\,,\ x<0\,,\ t>0\,.$$

Considering the set

$$\Omega_L = \{\lambda : \operatorname{Im}\lambda\leq 0\ \&\ \operatorname{Re}\omega_L(\lambda) = \operatorname{Re}(\alpha_L\lambda^2+\beta_L\lambda^4)\leq 0\}\ \text{(see fig 7)}$$

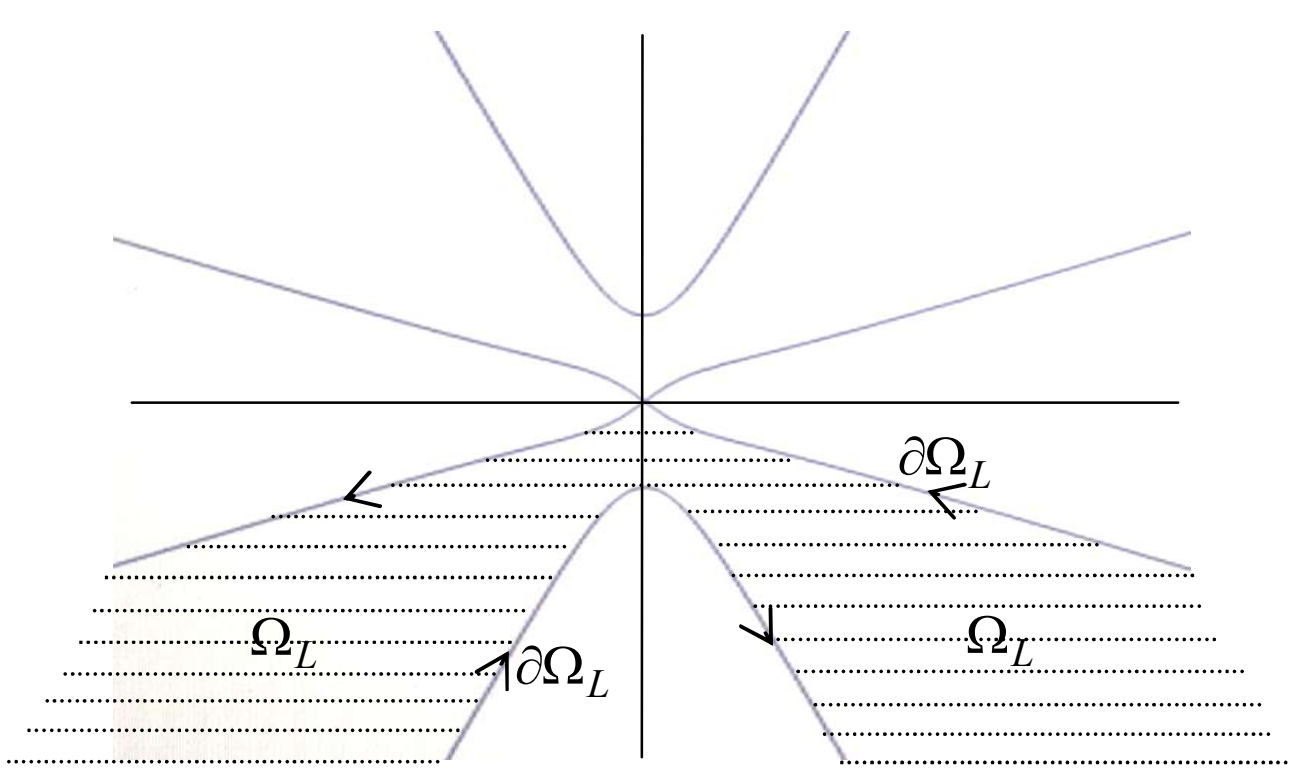

**Fig 7** *The set* $\Omega_L$ *and its boundary* $\partial\Omega_L$.

we write (3.1) as follows:

(3.3) $$2\pi u^L(x,t)=\int_{-\infty}^{\infty}e^{i\lambda x-\omega_L(\lambda)t}\hat{u}_0^L(\lambda)d\lambda+i\int_{\partial\Omega_L}\lambda e^{i\lambda x-\omega_L(\lambda)t}(\beta_L\lambda^2+\alpha_L)\widetilde{g}_0(\omega_L(\lambda),t)d\lambda$$
$$+\int_{\partial\Omega_L}e^{i\lambda x-\omega_L(\lambda)t}(\beta_L\lambda^2+\alpha_L)\widetilde{g}_{1,L}(\omega_L(\lambda),t)d\lambda-\beta_L i\int_{\partial\Omega_L}e^{i\lambda x-\omega_L(\lambda)t}\lambda\widetilde{g}_{2,L}(\omega_L(\lambda),t)d\lambda$$
$$-\beta_L\int_{\partial\Omega_L}e^{i\lambda x-\omega_L(\lambda)t}\widetilde{g}_{3,L}(\omega_L(\lambda),t)d\lambda+\int_{-\infty}^{\infty}e^{i\lambda x-\omega_L(\lambda)t}\widetilde{\hat{f}}_L(\lambda,\omega_L(\lambda),t)d\lambda,\ x<0,\ t>0.$$

Next we consider the mappings $\mu=\psi_L(\lambda)$ and $\lambda=\varphi_L(\mu)$ so that
$$\omega_L(\lambda)=\alpha_L\lambda^2+\beta_L\lambda^4=\mu^4=\omega_0(\mu),$$
which are defined as follows:
$$\mu=\psi_L(\lambda)=\begin{cases}-i\sqrt[4]{\alpha_L\lambda^2+\beta_L\lambda^4}\ \textit{if}\ \operatorname{Re}\lambda>0\\ -\sqrt[4]{\alpha_L\lambda^2+\beta_L\lambda^4}\ \textit{if}\ \operatorname{Re}\lambda<0\end{cases}\quad\text{and}\quad \lambda=\varphi_L(\mu)=\begin{cases}-\sqrt{\dfrac{-\alpha_L-\sqrt{4\beta_L\mu^4+\alpha_L^2}}{2\beta_L}}\ \textit{if}\ \operatorname{Re}\mu>0\\ -\sqrt{\dfrac{-\alpha_L+\sqrt{4\beta_L\mu^4+\alpha_L^2}}{2\beta_L}}\ \textit{if}\ \operatorname{Re}\mu<0.\end{cases}$$

(The choices of the functions $\sqrt[4]{\cdots}$ and $\sqrt{\cdots}$ are as in fig 2 and fig 3.)

***Remark*** As in *Section 2*,
$$\lim_{\substack{\mu\to0\\ \arg\mu=9\pi/8\ or\ \arg\mu=15\pi/8}}\varphi_L(\mu)=0\quad,\quad\lim_{\substack{\mu\to0\\ \arg\mu=11\pi/8\ or\ \arg\mu=13\pi/8}}\varphi_L(\mu)=-i\sqrt{\alpha_L/\beta_L},$$
and

(3.4) $$\frac{d\varphi_L(\mu)}{d\mu}=\mathrm{O}(\mu^2),\text{ as }\mu\to0\text{ with }\mu\in\Delta_L\text{ and }\frac{d\varphi_L(\mu)}{d\mu}=\mathrm{O}(1),\text{ as }\mu\to\infty\text{ with }\mu\in\Delta_L.$$

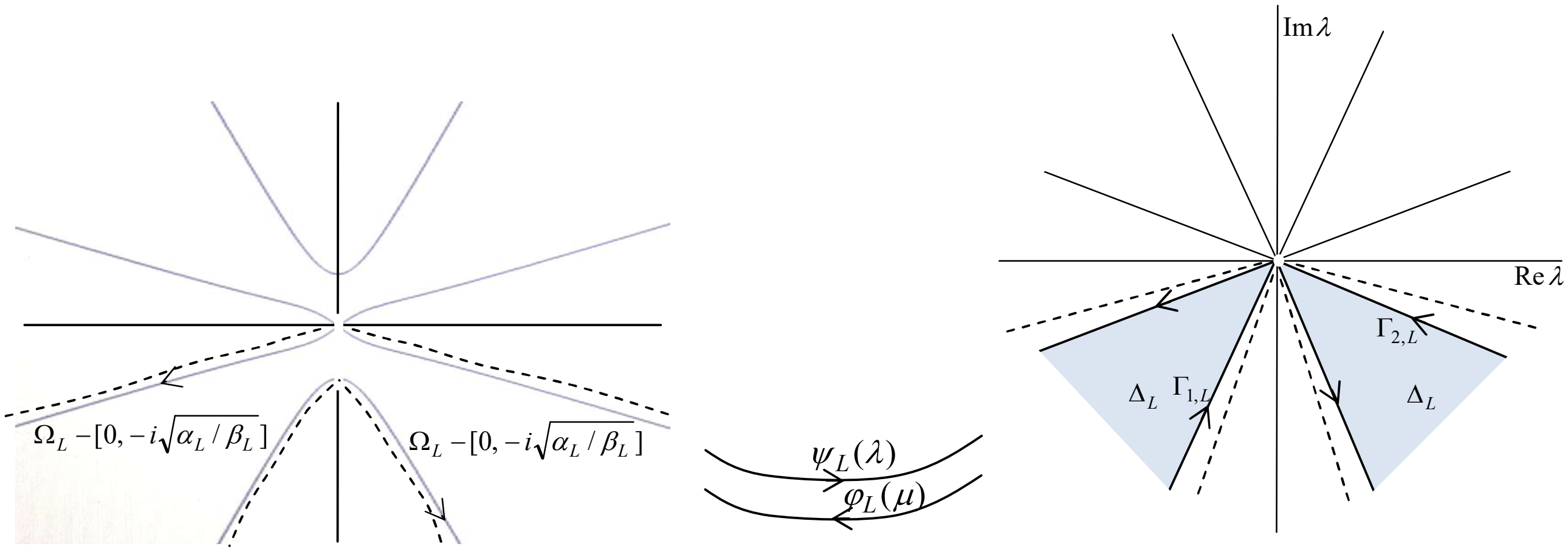

**Fig.8** *The mappings* $\mu=\psi_L(\lambda)$ *and* $\lambda=\varphi_L(\mu)$.

Then $\psi_L=\psi_L(\lambda)$ is a biholomorphic mapping from an open neighborhood of $\Omega_L-[0,-i\sqrt{\alpha_L/\beta_L}]$ to an open neighborhood of the set $\Delta_L:=\{\mu\in\mathbb{C}-\{0\}:\frac{9\pi}{8}\le\arg\mu\le\frac{11\pi}{8}\ or\ \frac{13\pi}{8}\le\arg\mu\le\frac{15\pi}{8}\}$, with inverse $\varphi_L=\varphi_L(\mu)$. (See fig 8.) Setting $\lambda=\varphi_R(\mu)$ and changing the variable from $\lambda$ to $\mu$, in the integrals over $\partial\Omega_L$, we write (3.3) in the following way:

(3.5) $$2\pi u^L(x,t)=\int_{-\infty}^{\infty}e^{i\lambda x-\omega_L(\lambda)t}\hat{u}_0^L(\lambda)d\lambda+i\int_{\Gamma_L}\varphi_L(\mu)e^{i\varphi_L(\mu)x-\omega_0(\mu)t}(\beta_L\varphi_L^2(\mu)+\alpha_L)\widetilde{g}_{0,L}(\omega_0(\mu),t)\frac{d\varphi_L(\mu)}{d\mu}d\mu$$

$$+\int_{\Gamma_L} e^{i\varphi_L(\mu)x-\omega_0(\mu)t}(\beta_L\varphi_L^2(\mu)+\alpha_L)\widetilde{g}_{1,L}(\omega_0(\mu),t)\frac{d\varphi_L(\mu)}{d\mu}d\mu-\beta_L i\int_{\Gamma_L}\varphi_L(\mu)e^{i\varphi_L(\mu)x-\omega_L(\lambda)t}\widetilde{g}_{2,L}(\omega_0(\mu),t)\frac{d\varphi_L(\mu)}{d\mu}d\mu$$

$$-\beta_L\int_{\Gamma_L} e^{i\varphi_L(\mu)x-\omega_0(\mu)t}\widetilde{g}_{3,L}(\omega_0(\mu),t)\frac{d\varphi_L(\mu)}{d\mu}d\mu+\int_{-\infty}^{\infty}e^{i\lambda x-\omega_L(\lambda)t}\widetilde{\widetilde{f}}_L(\lambda,\omega_L(\lambda),t)d\lambda\ ,\ x<0\ ,\ t>0\ .$$

where $\Gamma_L := \partial\Delta_L = \Gamma_{1,L}+\Gamma_{2,L}$ with $\Gamma_{1,L} := \{\lambda : \arg\lambda = \frac{9\pi}{8}\ or\ \frac{15\pi}{8}\}$ and $\Gamma_{2,L} := \{\lambda : \arg\lambda = \frac{11\pi}{8}\ or\ \frac{13\pi}{8}\}$. (See fig 9.)

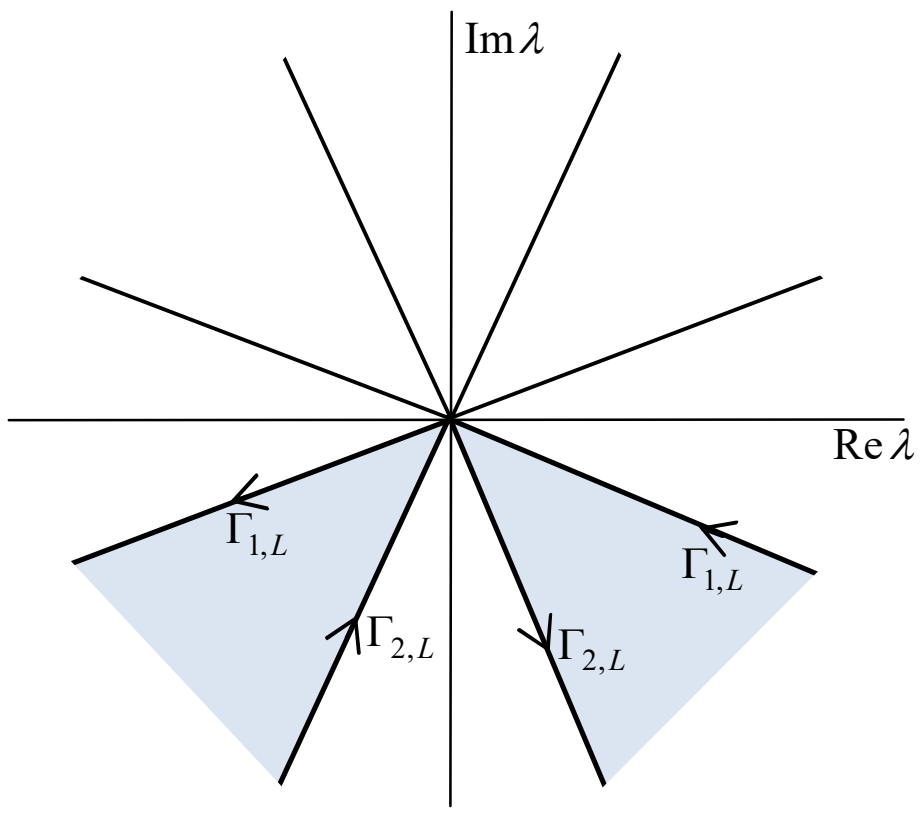

**Fig.9** *The contours* $\Gamma_{1,L}$ *and* $\Gamma_{2,L}$*, and their orientation.*

Next we write (3.5) in the following way:

$$(3.6)\quad 2\pi u^L(x,t)=\int_{-\infty}^{\infty}e^{i\lambda x-\omega_L(\lambda)t}\hat{u}_0^L(\lambda)d\lambda+i\int_{\Gamma_R}\varphi_L(-\mu)e^{i\varphi_L(-\mu)x-\omega_0(\mu)t}(\beta_L\varphi_L^2(-\mu)+\alpha_L)\widetilde{g}_0(\omega_0(\mu),t)\left\{-\frac{d[\varphi_L(-\mu)]}{d\mu}\right\}d\mu$$

$$+\int_{\Gamma_R} e^{i\varphi_L(-\mu)x-\omega_0(\mu)t}(\beta_L\varphi_L^2(-\mu)+\alpha_L)\widetilde{g}_{1,L}(\omega_0(\mu),t)\left\{-\frac{d[\varphi_L(-\mu)]}{d\mu}\right\}d\mu$$

$$-\beta_L i\int_{\Gamma_R}\varphi_L(-\mu)e^{i\varphi_L(-\mu)x-\omega_0(\mu)t}\widetilde{g}_{2,L}(\omega_0(\mu),t)\left\{-\frac{d[\varphi_L(-\mu)]}{d\mu}\right\}d\mu$$

$$-\beta_L\int_{\Gamma_R} e^{i\varphi_L(-\mu)x-\omega_0(\mu)t}\widetilde{g}_{3,L}(\omega_0(\mu),t)\left\{-\frac{d[\varphi_L(-\mu)]}{d\mu}\right\}d\mu+\int_{-\infty}^{\infty}e^{i\varphi_L(\mu)x-\omega_0(\mu)t}\widetilde{\widetilde{f}}_L(\lambda,\omega_0(\mu),t)d\lambda\ ,\ x<0\ ,\ t>0\ .$$

Also (3.1) gives

$$(3.7)\quad \hat{u}^L(-\varphi_L(-\mu),t)e^{\omega_0(\mu)t}=\hat{u}_0^L(-\varphi_L(-\mu))$$
$$-i\varphi_L(-\mu)(\beta_L\varphi_L^2(-\mu)+\alpha_L)\widetilde{g}_{0,L}(\omega_0(\mu),t)+(\beta_L\varphi_L^2(-\mu)+\alpha_L)\widetilde{g}_{1,L}(\omega_0(\mu),t)$$
$$+\beta_L i\varphi_L(-\mu)\widetilde{g}_{2,L}(\omega_0(\mu),t)-\beta_L\widetilde{g}_{3,L}(\omega_0(\mu),t)+\widetilde{\hat{f}}_L(-\varphi_L(-\mu),\omega_0(\mu),t)\ ,\ \mu\in\Delta_R\ .$$

$$(3.8)\quad \hat{u}^L(i\sigma_L(-\mu),t)e^{\omega_0(\mu)t}=\hat{u}_0^L(i\sigma_L(-\mu))+\beta_L\varphi_L^2(-\mu)\sigma_L(-\mu)\widetilde{g}_{0,L}(\omega_0(\mu),t)-\beta_L\varphi_L^2(-\mu)\widetilde{g}_{1,L}(\omega_0(\mu),t)$$
$$+\beta_L\sigma_L(-\mu)\widetilde{g}_{2,L}(\omega_0(\mu),t)-\beta_L\widetilde{g}_{3,L}(\omega_0(\mu),t)+\widetilde{\hat{f}}_L(i\sigma_L(-\mu),\omega_0(\mu),t)\ ,\ \mu\in\Delta_R\ .$$

where

$$\rho_L(\lambda)=\sqrt{\lambda^2+\frac{\alpha_L}{\beta_L}}\ \text{ and }\ \sigma_L(-\mu):=\rho_L(\varphi_L(-\mu))\ .$$

The function $\rho_L(\lambda)$ is defined as $\rho_R(\lambda)$ was defined in *Section 2* and it is analytic in

$$\Theta_L := \mathbb{C} - \{(-i\infty, -\sqrt{\tfrac{\alpha_L}{\beta_L}}i] \cup [\sqrt{\tfrac{\alpha_L}{\beta_L}}i, i\infty)\}.$$

## 4. Computing the integrals containing the unknown quantities

It follows from (1.3) that, for $\mu \in \mathbb{C}$ and $t > 0$,

$$(4.1)\qquad \begin{cases} \gamma_{11}\widetilde{g}_{0,R} + \gamma_{12}\widetilde{g}_{1,R} + \gamma_{13}\widetilde{g}_{2,R} + \gamma_{14}\widetilde{g}_{3,R} + \gamma_{15}\widetilde{g}_{0,L} + \gamma_{16}\widetilde{g}_{0,L} + \gamma_{17}\widetilde{g}_{0,L} + \gamma_{18}\widetilde{g}_{0,L} = \widetilde{I}_1(\omega_0(\mu),t) \\ \gamma_{21}\widetilde{g}_{0,R} + \gamma_{22}\widetilde{g}_{1,R} + \gamma_{23}\widetilde{g}_{2,R} + \gamma_{24}\widetilde{g}_{3,R} + \gamma_{25}\widetilde{g}_{0,L} + \gamma_{26}\widetilde{g}_{0,L} + \gamma_{27}\widetilde{g}_{0,L} + \gamma_{28}\widetilde{g}_{0,L} = \widetilde{I}_2(\omega_0(\mu),t) \\ \gamma_{31}\widetilde{g}_{0,R} + \gamma_{32}\widetilde{g}_{1,R} + \gamma_{33}\widetilde{g}_{2,R} + \gamma_{34}\widetilde{g}_{3,R} + \gamma_{35}\widetilde{g}_{0,L} + \gamma_{36}\widetilde{g}_{0,L} + \gamma_{37}\widetilde{g}_{0,L} + \gamma_{38}\widetilde{g}_{0,L} = \widetilde{I}_3(\omega_0(\mu),t) \\ \gamma_{41}\widetilde{g}_{0,R} + \gamma_{42}\widetilde{g}_{1,R} + \gamma_{43}\widetilde{g}_{2,R} + \gamma_{44}\widetilde{g}_{3,R} + \gamma_{45}\widetilde{g}_{0,L} + \gamma_{46}\widetilde{g}_{0,L} + \gamma_{47}\widetilde{g}_{0,L} + \gamma_{48}\widetilde{g}_{0,L} = \widetilde{I}_4(\omega_0(\mu),t), \end{cases}$$

where we used the abbreviations $\widetilde{g}_{j,R} = \widetilde{g}_{j,R}(\omega_0(\mu),t)$ and $\widetilde{g}_{j,L} = \widetilde{g}_{j,L}(\omega_0(\mu),t)$.

Now we proceed to solve the system of equations (2.6), (2.7), (3.7), (3.8) and (4.1), in the unknown quantities

$$(4.2)\qquad \widetilde{g}_{0,R}(\omega_0(\mu),t),\ \widetilde{g}_{1,R}(\omega_0(\mu),t),\ \widetilde{g}_{2,R}(\omega_0(\mu),t),\ \widetilde{g}_{3,R}(\omega_0(\mu),t),$$

and

$$(4.3)\qquad \widetilde{g}_{0,L}(\omega_0(\mu),t),\ \widetilde{g}_{1,L}(\omega_0(\mu),t),\ \widetilde{g}_{2,L}(\omega_0(\mu),t),\ \widetilde{g}_{3,L}(\omega_0(\mu),t),$$

with $\mu \in \Delta_R$ and $t > 0$.

The matrix of the coefficients of the unknowns is

$$\mathfrak{D}(\mu) := \begin{bmatrix} i\varphi_R(\mu)\Phi_R(\mu) & -\Phi_R(\mu) & -\beta_R i\varphi_R(\mu) & \beta_R & 0 & 0 & 0 & 0 \\ -\beta_R\varphi_R^2(\mu)\sigma_R(\mu) & \beta_R\varphi_R^2(\mu) & -\beta_R\sigma_R(\mu) & \beta_R & 0 & 0 & 0 & 0 \\ 0 & 0 & 0 & 0 & -i\varphi_L(-\mu)\Phi_L(-\mu) & \Phi_L(-\mu) & \beta_L i\varphi_L(-\mu) & -\beta_L \\ 0 & 0 & 0 & 0 & \beta_L\varphi_L^2(-\mu)\sigma_L(-\mu) & -\beta_L\varphi_L^2(-\mu) & \beta_L\sigma_L(-\mu) & -\beta_L \\ \gamma_{11} & \gamma_{12} & \gamma_{13} & \gamma_{14} & \gamma_{15} & \gamma_{16} & \gamma_{17} & \gamma_{18} \\ \gamma_{21} & \gamma_{22} & \gamma_{23} & \gamma_{24} & \gamma_{25} & \gamma_{26} & \gamma_{27} & \gamma_{28} \\ \gamma_{31} & \gamma_{32} & \gamma_{33} & \gamma_{34} & \gamma_{35} & \gamma_{36} & \gamma_{37} & \gamma_{38} \\ \gamma_{41} & \gamma_{42} & \gamma_{43} & \gamma_{44} & \gamma_{45} & \gamma_{46} & \gamma_{47} & \gamma_{48} \end{bmatrix},$$

where $\Phi_R(\mu) := \beta_R\varphi_R^2(\mu) + \alpha_R$ and $\Phi_L(\mu) := \beta_L\varphi_L^2(\mu) + \alpha_L$.

***Remark*** A calculation shows that

$$\lim_{\substack{\mu\to 0 \\ \arg\mu=\pi/8 \text{ or } \arg\mu=7\pi/8}} \det\mathfrak{D}(\mu) = \begin{bmatrix} 0 & -\alpha_R & 0 & \beta_R & 0 & 0 & 0 & 0 \\ 0 & 0 & -\beta_R^{1/2}\alpha_R^{1/2} & \beta_R & 0 & 0 & 0 & 0 \\ 0 & 0 & 0 & 0 & 0 & \alpha_L & 0 & -\beta_L \\ 0 & 0 & 0 & 0 & 0 & 0 & \alpha_L^{1/2}\beta_L^{1/2} & -\beta_L \\ \gamma_{11} & \gamma_{12} & \gamma_{13} & \gamma_{14} & \gamma_{15} & \gamma_{16} & \gamma_{17} & \gamma_{18} \\ \gamma_{21} & \gamma_{22} & \gamma_{23} & \gamma_{24} & \gamma_{25} & \gamma_{26} & \gamma_{27} & \gamma_{28} \\ \gamma_{31} & \gamma_{32} & \gamma_{33} & \gamma_{34} & \gamma_{35} & \gamma_{36} & \gamma_{37} & \gamma_{38} \\ \gamma_{41} & \gamma_{42} & \gamma_{43} & \gamma_{44} & \gamma_{45} & \gamma_{46} & \gamma_{47} & \gamma_{48} \end{bmatrix}$$

and

$$\lim_{\substack{\mu\to 0\\ \arg\mu=3\pi/8 \text{ or } \arg\mu=5\pi/8}} \det\mathfrak{D}(\mu) = \begin{bmatrix} 0 & 0 & \beta_R^{1/2}\alpha_R^{1/2} & \beta_R & 0 & 0 & 0 & 0 \\ 0 & -\beta_R^{1/2}\alpha_R^{1/2} & 0 & \beta_R & 0 & 0 & 0 & 0 \\ 0 & 0 & 0 & 0 & 0 & 0 & \alpha_L^{1/2}\beta_L^{1/2} & -\beta_L \\ 0 & 0 & 0 & 0 & 0 & \alpha_L & 0 & -\beta_L \\ \gamma_{11} & \gamma_{12} & \gamma_{13} & \gamma_{14} & \gamma_{15} & \gamma_{16} & \gamma_{17} & \gamma_{18} \\ \gamma_{21} & \gamma_{22} & \gamma_{23} & \gamma_{24} & \gamma_{25} & \gamma_{26} & \gamma_{27} & \gamma_{28} \\ \gamma_{31} & \gamma_{32} & \gamma_{33} & \gamma_{34} & \gamma_{35} & \gamma_{36} & \gamma_{37} & \gamma_{38} \\ \gamma_{41} & \gamma_{42} & \gamma_{43} & \gamma_{44} & \gamma_{45} & \gamma_{46} & \gamma_{47} & \gamma_{48} \end{bmatrix}.$$

In particular, if either of the above limits is not zero, the function $\det\mathfrak{D}(\mu)$ is not identically zero.

Setting

$$J_{1,R}^*(\mu,t) := \hat{u}^R(-\varphi_R(\mu),t)e^{\omega_0(\mu)t}, \; J_{1,R}(\mu,t) := -\hat{u}_0^R(-\varphi_R(\mu)) - \tilde{\hat{f}}_R(-\varphi_R(\mu),\omega_0(\mu),t)$$
$$J_{2,R}^*(\mu,t) := \hat{u}_R(i\sigma_R(\mu),t)e^{\omega_0(\mu)t}, \; J_{2,R}(\mu,t) := -\hat{u}_0^R(i\sigma_R(\mu)) - \tilde{\hat{f}}_R(i\sigma_R(\mu),\omega_0(\mu),t)$$
$$J_{1,L}^*(\mu,t) := \hat{u}^L(-\varphi_L(-\mu),t)e^{\omega_0(\mu)t}, \; J_{1,L}(\mu,t) := -\hat{u}_0^L(-\varphi_L(-\mu)) - \tilde{\hat{f}}_L(-\varphi_L(-\mu),\omega_0(\mu),t)$$
$$J_{2,L}^*(\mu,t) := \hat{u}^L(i\sigma_L(-\mu),t)e^{\omega_0(\mu)t}, \; J_{2,L}(\mu,t) := -\hat{u}_0^L(i\sigma_L(-\mu)) - \tilde{\hat{f}}_L(i\sigma_L(-\mu),\omega_0(\mu),t),$$

we write the Cramer solutions for the unknowns (4.2) and (4.3) of the system of equations {(2.6), (2.7), (3.7), (3.8), (4.1)} as follows: For $\mu\in\Delta_R$,

$$G_{j,R}^*(\mu,t) := \tilde{g}_{j,R}(\omega_0(\mu),t) = \frac{\det\mathfrak{D}_{j,R}^*(\mu,t)}{\det\mathfrak{D}(\mu)}, \; j=0,1,2,3,$$

and

$$G_{j,L}^*(\mu,t) := \tilde{g}_{j,L}(\omega_0(\mu),t) = \frac{\det\mathfrak{D}_{j,L}^*(\mu,t)}{\det\mathfrak{D}(\mu)}, \; j=0,1,2,3,$$

where $\mathfrak{D}_{j,R}^*(\mu,t)$ is the matrix which results from $\mathfrak{D}(\mu)$, by replacing its $(j+1)^{th}-column$ ( $j=0,1,2,3$ ) by

$$\mathfrak{C}^*(\mu,t) := \begin{bmatrix} J_{1,R}(\mu,t)+J_{1,R}^*(\mu,t) \\ J_{2,R}(\mu,t)+J_{2,R}^*(\mu,t) \\ J_{1,L}(\mu,t)+J_{1,L}^*(\mu,t) \\ J_{2,L}(\mu,t)+J_{2,L}^*(\mu,t) \\ \tilde{I}_1(\omega_0(\mu),t) \\ \tilde{I}_2(\omega_0(\mu),t) \\ \tilde{I}_3(\omega_0(\mu),t) \\ \tilde{I}_4(\omega_0(\mu),t) \end{bmatrix}$$

and $\mathfrak{D}_{j,L}^*(\mu,t)$ is the matrix which results again from $\mathfrak{D}(\mu)$, by replacing its $(j+5)^{th}-column$ ( $j=0,1,2,3$ ) by the column $\mathfrak{C}^*(\mu,t)$.

For example,

$$\mathfrak{D}_{0,R}^{*}(\mu)=\begin{bmatrix} J_{1,R}(\mu,t)+J_{1,R}^{*}(\mu,t) & -\Phi_R(\mu) & -\beta_R i\varphi_R(\mu) & \beta_R & 0 & 0 & 0 & 0 \\ J_{2,R}(\mu,t)+J_{2,R}^{*}(\mu,t) & \beta_R\varphi_R^2(\mu) & -\beta_R\sigma_R(\mu) & \beta_R & 0 & 0 & 0 & 0 \\ J_{1,L}(\mu,t)+J_{1,L}^{*}(\mu,t) & 0 & 0 & 0 & -i\varphi_L(-\mu)\Phi_L(-\mu) & \Phi_L(-\mu) & \beta_L i\varphi_L(-\mu) & -\beta_L \\ J_{2,L}(\mu,t)+J_{2,L}^{*}(\mu,t) & 0 & 0 & 0 & \beta_L\varphi_L^2(-\mu)\sigma_L(-\mu) & -\beta_L\varphi_L^2(-\mu) & \beta_L\sigma_L(-\mu) & -\beta_L \\ \widetilde{I}_1(\omega_0(\mu),t) & \gamma_{12} & \gamma_{13} & \gamma_{14} & \gamma_{15} & \gamma_{16} & \gamma_{17} & \gamma_{18} \\ \widetilde{I}_2(\omega_0(\mu),t) & \gamma_{22} & \gamma_{23} & \gamma_{24} & \gamma_{25} & \gamma_{26} & \gamma_{27} & \gamma_{28} \\ \widetilde{I}_3(\omega_0(\mu),t) & \gamma_{32} & \gamma_{33} & \gamma_{34} & \gamma_{35} & \gamma_{36} & \gamma_{37} & \gamma_{38} \\ \widetilde{I}_4(\omega_0(\mu),t) & \gamma_{42} & \gamma_{43} & \gamma_{44} & \gamma_{45} & \gamma_{46} & \gamma_{47} & \gamma_{48} \end{bmatrix}$$

and

$$\mathfrak{D}_{2,L}^{*}(\mu)=\begin{bmatrix} i\varphi_R(\mu)\Phi_R(\mu) & J_{1,R}(\mu,t)+J_{1,R}^{*}(\mu,t) & -\beta_R i\varphi_R(\mu) & \beta_R & 0 & 0 & 0 & 0 \\ -\beta_R\varphi_R^2(\mu)\sigma_R(\mu) & J_{2,R}(\mu,t)+J_{2,R}^{*}(\mu,t) & -\beta_R\sigma_R(\mu) & \beta_R & 0 & 0 & 0 & 0 \\ 0 & J_{1,L}(\mu,t)+J_{1,L}^{*}(\mu,t) & 0 & 0 & -i\varphi_L(-\mu)\Phi_L(-\mu) & \Phi_L(-\mu) & \beta_L i\varphi_L(-\mu) & -\beta_L \\ 0 & J_{2,L}(\mu,t)+J_{2,L}^{*}(\mu,t) & 0 & 0 & \beta_L\varphi_L^2(-\mu)\sigma_L(-\mu) & -\beta_L\varphi_L^2(-\mu) & \beta_L\sigma_L(-\mu) & -\beta_L \\ \gamma_{11} & \widetilde{I}_1(\omega_0(\mu),t) & \gamma_{13} & \gamma_{14} & \gamma_{15} & \gamma_{16} & \gamma_{17} & \gamma_{18} \\ \gamma_{21} & \widetilde{I}_2(\omega_0(\mu),t) & \gamma_{23} & \gamma_{24} & \gamma_{25} & \gamma_{26} & \gamma_{27} & \gamma_{28} \\ \gamma_{31} & \widetilde{I}_3(\omega_0(\mu),t) & \gamma_{33} & \gamma_{34} & \gamma_{35} & \gamma_{36} & \gamma_{37} & \gamma_{38} \\ \gamma_{41} & \widetilde{I}_4(\omega_0(\mu),t) & \gamma_{43} & \gamma_{44} & \gamma_{45} & \gamma_{46} & \gamma_{47} & \gamma_{48} \end{bmatrix}$$

Finally, for $\mu\in\Delta_R$, we define

$$G_{j,R}(\mu,t):=\frac{\det\mathfrak{D}_{j,R}(\mu,t)}{\det\mathfrak{D}(\mu)},\ G_{j,L}(\mu,t):=\frac{\det\mathfrak{D}_{j,L}(\mu,t)}{\det\mathfrak{D}(\mu)},\ j=0,1,2,3,$$

where $\mathfrak{D}_{j,R}(\mu,t)$ is the matrix which results from $\mathfrak{D}(\mu)$, by replacing its $(j+1)^{th}-column$ by

$$\mathfrak{C}(\mu,t):=\begin{bmatrix} J_{1,R}(\mu,t) \\ J_{2,R}(\mu,t) \\ J_{1,L}(\mu,t) \\ J_{2,L}(\mu,t) \\ \widetilde{I}_1(\omega_0(\mu),t) \\ \widetilde{I}_2(\omega_0(\mu),t) \\ \widetilde{I}_3(\omega_0(\mu),t) \\ \widetilde{I}_4(\omega_0(\mu),t) \end{bmatrix}$$

and $\mathfrak{D}_{j,L}(\mu,t)$ is the matrix which results from $\mathfrak{D}(\mu)$, by replacing its $(j+5)^{th}-column$ by $\mathfrak{C}(\mu,t)$, $j=0,1,2,3$.

***Assumption which guarantees that*** $\det\mathfrak{D}(\mu)\neq 0$ ***for large*** $|\mu|$**.** Let us observe that

$$w_1:=\varphi_R(\mu)\ \Rightarrow\ (2\beta_R w_1^2-\alpha_R)^2=4\beta_R\mu^4+\alpha_R^2,$$

$$w_2:=\varphi_L(-\mu)\ \Rightarrow\ (2\beta_L w_2^2-\alpha_L)^2=4\beta_L\mu^4+\alpha_L^2,$$

$$w_3:=\sigma_R(\mu)=\rho_R(\varphi_R(\mu))\ \Rightarrow\ w_3^2=w_1^2+\frac{\alpha_R}{\beta_R},$$

$$w_4:=\sigma_L(-\mu)=\rho_L(\varphi_R(-\mu))\ \Rightarrow\ w_4^2=w_2^2+\frac{\alpha_L}{\beta_L}.$$

Thus, considering the polynomials $p_\ell = p_\ell(w_1, w_2, w_3, w_4, \mu)$, $\ell = 1, 2, 3, 4, 5$, in the variables $w_1, w_2, w_3, w_4, \mu$, with

$$p_1 := (2\beta_R w_1^2 - \alpha_R)^2 - 4\beta_R \mu^4 - \alpha_R^2,$$
$$p_2 := (2\beta_L w_2^2 - \alpha_L)^2 - 4\beta_L \mu^4 - \alpha_L^2,$$
$$p_3 := w_3^2 - w_1^2 - \frac{\alpha_R}{\beta_R},$$
$$p_4 := w_4^2 - w_2^2 - \frac{\alpha_L}{\beta_L}$$

and

$$p_5 := \det \mathfrak{D}(\mu)$$

$$= \det \begin{bmatrix} iw_1(\beta_R w_1^2 + \alpha_R) & -(\beta_R w_1^2 + \alpha_R) & -\beta_R i w_1 & \beta_R & 0 & 0 & 0 & 0 \\ -\beta_R w_1^2 w_3 & \beta_R w_1^2 & -\beta_R w_3 & \beta_R & 0 & 0 & 0 & 0 \\ 0 & 0 & 0 & 0 & -iw_2(\beta_L w_2^2 + \alpha_L) & \beta_L w_2^2 + \alpha_L & \beta_L i w_2 & -\beta_L \\ 0 & 0 & 0 & 0 & \beta_L w_2^2 w_4^2 & -\beta_L w_2^2 & \beta_L w_4 & -\beta_L \\ \gamma_{11} & \gamma_{12} & \gamma_{13} & \gamma_{14} & \gamma_{15} & \gamma_{16} & \gamma_{17} & \gamma_{18} \\ \gamma_{21} & \gamma_{22} & \gamma_{23} & \gamma_{24} & \gamma_{25} & \gamma_{26} & \gamma_{27} & \gamma_{28} \\ \gamma_{31} & \gamma_{32} & \gamma_{33} & \gamma_{34} & \gamma_{35} & \gamma_{36} & \gamma_{37} & \gamma_{38} \\ \gamma_{41} & \gamma_{42} & \gamma_{43} & \gamma_{44} & \gamma_{45} & \gamma_{46} & \gamma_{47} & \gamma_{48} \end{bmatrix},$$

we see that if the complex number $\mu_0$ satisfies the equation $\det \mathfrak{D}(\mu_0) = 0$, then the corresponding numbers $w_{1,0}, w_{2,0}, w_{3,0}, w_{4,0}, \mu_0$, satisfy the system of algebraic equations: $p_\ell(w_1, w_2, w_3, w_4, \mu) = 0$, $\ell = 1, 2, 3, 4, 5$. Eliminating $w_1, w_2, w_3, w_4$ between these equations, we obtain a polynomial $P(\mu)$ such that $P(\mu_0) = 0$. It follows that if the polynomial $P(\mu)$ is not identically zero then the function $\det \mathfrak{D}(\mu)$ has finitely many zeros, and, therefore, there exists an upper bound for the absolute values of these zeros. Thus, if $P(\mu) \not\equiv 0$, there exists a positive number $\mathrm{A}$ such that

(4.4) $$|\mu| \geq \mathrm{A} \text{ with } \mu \in \Gamma_R \Rightarrow \det \mathfrak{D}(\mu) \neq 0.$$

Also, since the polynomial $P(\mu)$ can be explicitly computed, an estimate of a number $\mathrm{A}$, with the above property, can also be found.

***Further deformation of the contours*** Assuming that $P(\mu) \not\equiv 0$ and choosing $\mathrm{A}$ so that (4.4) holds, we consider the contour

$$\mathrm{E_A} := \{\mu : |\mu| \geq \mathrm{A}, \arg \mu = \tfrac{7\pi}{8}\} + \{\mu : |\mu| = \mathrm{A}, \tfrac{7\pi}{8} \geq \arg \mu \geq \tfrac{5\pi}{8}\} + \{\mu : |\mu| \geq \mathrm{A}, \arg \mu = \tfrac{5\pi}{8}\}$$
$$+ \{\mu : |\mu| \geq \mathrm{A}, \arg \mu = \tfrac{3\pi}{8}\} + \{\mu : |\mu| = \mathrm{A}, \tfrac{3\pi}{8} \geq \arg \mu \geq \tfrac{\pi}{8}\} + \{\mu : |\mu| \geq \mathrm{A}, \arg \mu = \tfrac{\pi}{8}\}.$$ (See fig 10.)

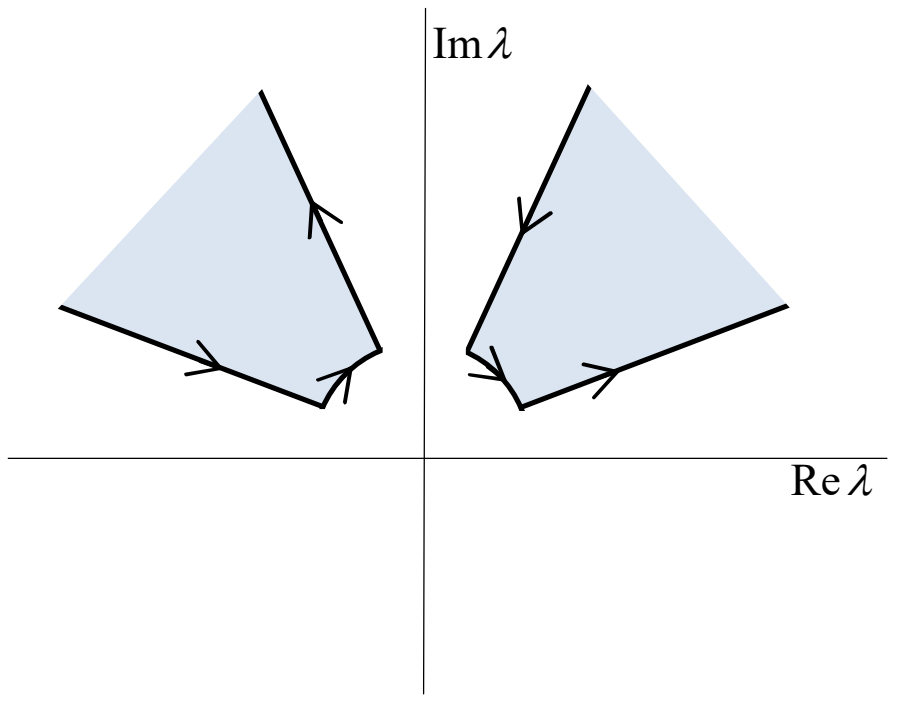

**Fig.10** *The parts of the contour* $\mathrm{E_A}$ *and their orientations.*

Thus, based on the above analysis, we can formulate the following:

**Theorem** *Suppose that the data satisfy (1.4), (1.5), (1.6) and that the parameters are such that the polynomial* $P(\mu)\not\equiv 0$. *Then the UTM-based solution to the interface problem {(1.1), (1.2) &(1.3)} is given by the following formulas:*

$$\text{(4.5)}\quad 2\pi u^R(x,t)=\int_{-\infty}^{\infty}e^{i\lambda x-\omega_R(\lambda)t}\hat{u}_0^R(\lambda)d\lambda-i\int_{\partial\mathcal{Z}_A^+}\lambda e^{i\lambda x-\omega_R(\lambda)t}(\beta_R\lambda^2+\alpha_R)G_{0,R}(\psi_R(\lambda),t)d\lambda$$
$$-\int_{\partial\mathcal{Z}_A^+}e^{i\lambda x-\omega_R(\lambda)t}(\beta_R\lambda^2+\alpha_R)G_{1,R}(\psi_R(\lambda),t)d\lambda+\beta_R i\int_{\partial\mathcal{Z}_A^+}\lambda e^{i\lambda x-\omega_R(\lambda)t}G_{2,R}(\psi_R(\lambda),t)d\lambda$$
$$+\beta_R\int_{\partial\mathcal{Z}_A^+}e^{i\lambda x-\omega_R(\lambda)t}G_{3,R}(\psi_R(\lambda),t)d\lambda+\int_{-\infty}^{\infty}e^{i\lambda x-\omega_R(\lambda)t}\tilde{\hat{f}}_R(\lambda,\omega_R(\lambda),t)d\lambda,\ x>0,\ t>0,$$

*and*

$$\text{(4.6)}\quad 2\pi u^L(x,t)=\int_{-\infty}^{\infty}e^{i\lambda x-\omega_L(\lambda)t}\hat{u}_0^L(\lambda)d\lambda+i\int_{\partial\mathcal{Z}_A^-}\lambda e^{i\lambda x-\omega_L(\lambda)t}(\beta_L\lambda^2+\alpha_L)G_{0,L}(-\psi_L(\lambda),t)d\lambda$$
$$+\int_{\partial\mathcal{Z}_A^-}e^{i\lambda x-\omega_L(\lambda)t}(\beta_L\lambda^2+\alpha_L)G_{1,L}(-\psi_L(\lambda),t)d\lambda-\beta_L i\int_{\partial\mathcal{Z}_A^-}\lambda e^{i\lambda x-\omega_L(\lambda)t}G_{2,L}(-\psi_L(\lambda),t)d\lambda$$
$$-\beta_L\int_{\partial\mathcal{Z}_A^-}e^{i\lambda x-\omega_L(\lambda)t}G_{3,L}(-\psi_L(\lambda),t)d\lambda+\int_{-\infty}^{\infty}e^{i\lambda x-\omega_L(\lambda)t}\tilde{\hat{f}}_L(\lambda,\omega_L(\lambda),t)d\lambda,\ x<0,\ t>0,$$

*where* $\mathcal{Z}_A^+:=\{\lambda\in\Omega_R:|\lambda|\ge A\}$ (fig 11), $\mathcal{Z}_A^-:=\{\lambda\in\Omega_L:|\lambda|\ge A\}$ (fig 12), *and* A *is sufficiently large.*

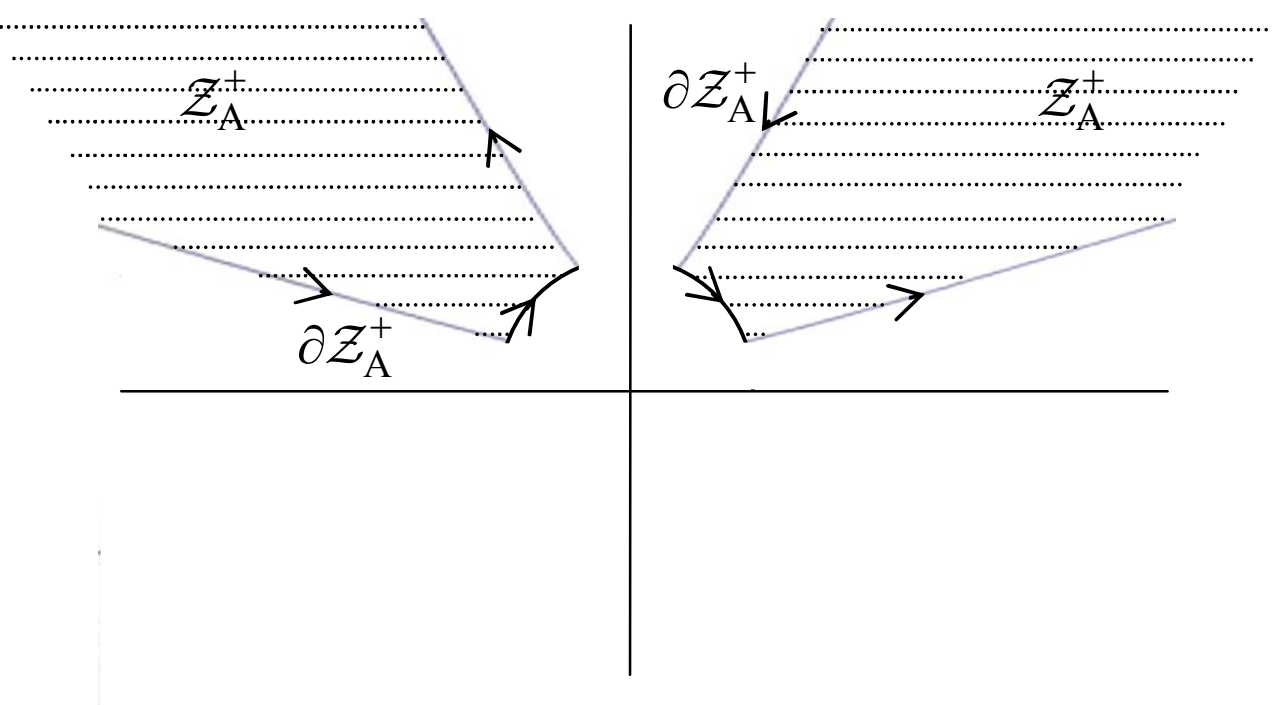

**Fig11** *The set* $\mathcal{Z}_A^+$ *and its boundary* $\partial\mathcal{Z}_A^+$.

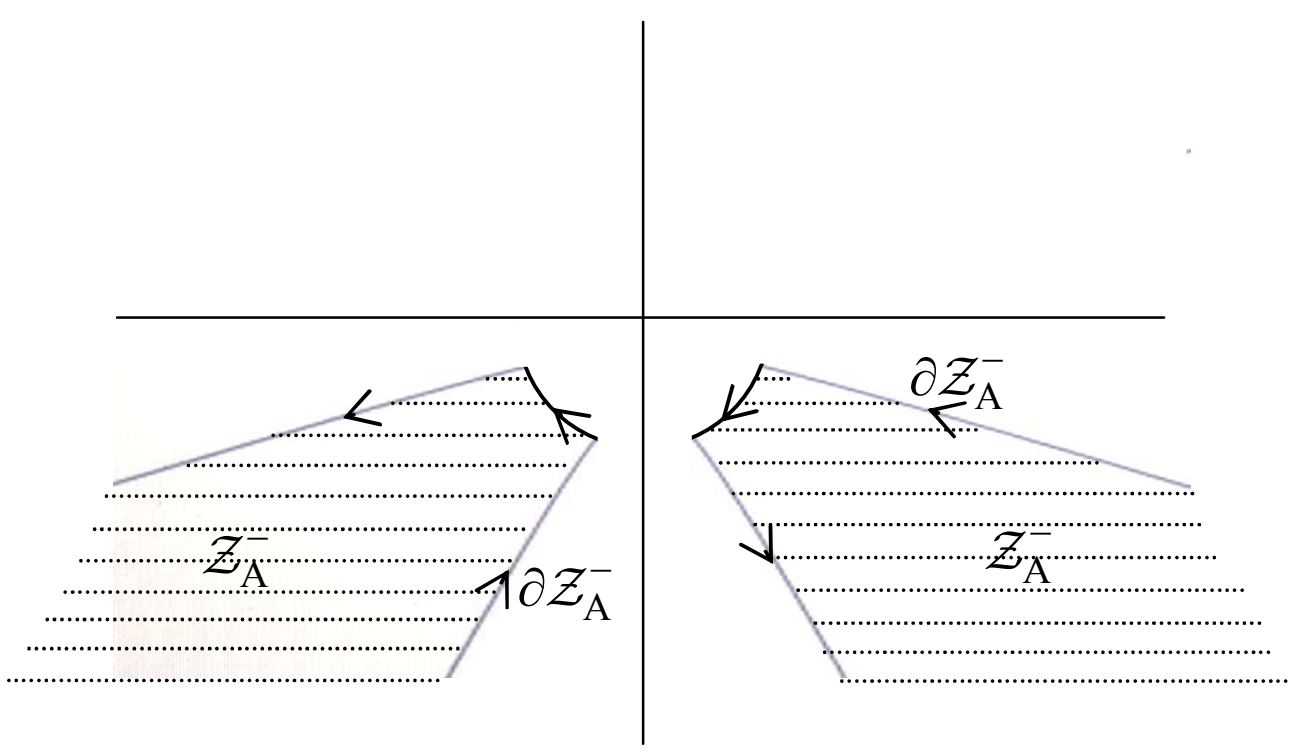

**Fig12** *The set* $\mathcal{Z}_A^-$ *and its boundary* $\partial\mathcal{Z}_A^-$.

***Proof*** Substituting $\widetilde{g}_{j,R}$ in (2.4) and deforming the contour $\Gamma_R$ to the contour $\mathrm{E}_\mathrm{A}$, in order to avoid the zeros of the function $\det\mathfrak{D}(\mu)$, we obtain the following formula:

$$(4.7)\quad 2\pi u^R(x,t) = \int_{-\infty}^{\infty} e^{i\lambda x-\omega_R(\lambda)t}\hat{u}_0^R(\lambda)d\lambda - i\int_{\mathrm{E}_\mathrm{A}} \varphi_R(\mu)e^{i\varphi_R(\mu)x-\omega_0(\mu)t}(\beta_R\varphi_R^2(\mu)+\alpha_R)G_{0,R}^*(\mu,t)\frac{d\varphi_R(\mu)}{d\mu}d\mu$$
$$-\int_{\mathrm{E}_\mathrm{A}} e^{i\varphi_R(\mu)x-\omega_0(\mu)t}(\beta_R\varphi_R^2(\mu)+\alpha_R)G_{1,R}^*(\mu,t)\frac{d\varphi_R(\mu)}{d\mu}d\mu$$
$$+\beta_R i\int_{\mathrm{E}_\mathrm{A}} \varphi_R(\mu)e^{i\varphi_R(\mu)x-\omega_0(\mu)t}G_{2,R}^*(\mu,t)\frac{d\varphi_R(\mu)}{d\mu}d\mu$$
$$+\beta_R\int_{\mathrm{E}_\mathrm{A}} e^{i\varphi_R(\mu)x-\omega_0(\mu)t}G_{3,R}^*(\mu,t)\frac{d\varphi_R(\mu)}{d\mu}d\mu + \int_{-\infty}^{\infty} e^{i\lambda x-\omega_R(\lambda)t}\widetilde{\widetilde{f}}_R(\lambda,\omega_R(\lambda),t)d\lambda\,,\ x>0\,,\ t>0,$$

where $\mathrm{A}$ is as in (4.4).
An analogous formula holds also for $u^L(x,t)$, in the case $x<0$, $t>0$, which is derived by substituting the solution $\widetilde{g}_{j,L}$ in (3.6).

Next, we observe that the coefficients of the quantities $J_{1,R}^*(\mu,t)$, $J_{2,R}^*(\mu,t)$, $J_{1,L}^*(\mu,t)$, $J_{2,L}^*(\mu,t)$, in the determinants $\det\mathfrak{D}_{j,R}^*(\mu)$ and $\det\mathfrak{D}_{j,L}^*(\mu)$, are $\mathrm{O}(\mu^7)$, as $\mu\to\infty$, with $\mu\in\Delta_R$. Also, keeping in mind that

$$\det\mathfrak{D}(\mu)=\mathrm{O}(\mu^{10})\,,\ \frac{d\varphi_R(\mu)}{d\mu}=\mathrm{O}(1)\,,\ \frac{d\varphi_L(-\mu)}{d\mu}=\mathrm{O}(1)\ \text{(see (2.5) and (3.4))},$$

we conclude that, as $\mu\to\infty$, with $\mu\in\Delta_R$,

$$(4.8)\quad e^{-\omega_0(\mu)t}G_{j,R}^*(\mu,t)-e^{-\omega_0(\mu)t}G_{j,R}(\mu,t)=\mathrm{O}(\mu^6),$$

and

$$(4.9)\quad e^{-\omega_0(\mu)t}G_{j,L}^*(\mu,t)-e^{-\omega_0(\mu)t}G_{j,L}(\mu,t)=\mathrm{O}(\mu^6).$$

Thus, setting

$$\mathcal{U}_{j,R}(\mu,t):=G_{j,R}^*(\mu,t)-G_{j,R}(\mu,t)\,,\ j=0,1,2,3\,,$$

we see that, in view of (4.8) and (4.9), and Cauchy's theorem,

$$\int_{\mathrm{E}_\mathrm{A}} \varphi_R(\mu)e^{i\varphi_R(\mu)x-\omega_0(\mu)t}(\beta_R\varphi_R^2(\mu)+\alpha_R)\mathcal{U}_{0,R}(\mu,t)\frac{d\varphi_R(\mu)}{d\mu}d\mu=0,$$
$$\int_{\mathrm{E}_\mathrm{A}} e^{i\varphi_R(\mu)x-\omega_0(\mu)t}(\beta_R\varphi_R^2(\mu)+\alpha_R)\mathcal{U}_{1,R}(\mu,t)\frac{d\varphi_R(\mu)}{d\mu}d\mu=0,$$
$$\int_{\mathrm{E}_\mathrm{A}} \varphi_R(\mu)e^{i\varphi_R(\mu)x-\omega_0(\mu)t}\mathcal{U}_{2,R}(\mu,t)\frac{d\varphi_R(\mu)}{d\mu}d\mu=0,$$
$$\int_{\mathrm{E}_\mathrm{A}} e^{i\varphi_R(\mu)x-\omega_0(\mu)t}\mathcal{U}_{3,R}(\mu,t)\frac{d\varphi_R(\mu)}{d\mu}d\mu=0.$$

Indeed, the above equations follow from the following form of Jordan's lemma:

$$\lim_{\mathrm{B}\to+\infty}\int_{\{\mu:\,|\mu|=\mathrm{B}\}\cap\Delta_R} e^{i\varphi_R(\mu)x}\mu^N d\mu=0\,,\ \text{for } x>0 \text{ and } N\in\mathbb{N}.$$

Therefore, (4.7) becomes

$$(4.10)\quad 2\pi u^R(x,t) = \int_{-\infty}^{\infty} e^{i\lambda x-\omega_R(\lambda)t}\hat{u}_0^R(\lambda)d\lambda - i\int_{\mathrm{E}_\mathrm{A}} \varphi_R(\mu)e^{i\varphi_R(\mu)x-\omega_0(\mu)t}(\beta_R\varphi_R^2(\mu)+\alpha_R)G_{0,R}(\mu,t)\frac{d\varphi_R(\mu)}{d\mu}d\mu$$

$$-\int_{E_A} e^{i\varphi_R(\mu)x-\omega_0(\mu)t}(\beta_R\varphi_R^2(\mu)+\alpha_R)G_{1,R}(\mu,t)\frac{d\varphi_R(\mu)}{d\mu}d\mu$$

$$+\beta_R i\int_{E_A}\varphi_R(\mu)e^{i\varphi_R(\mu)x-\omega_0(\mu)t}G_{2,R}(\mu,t)\frac{d\varphi_R(\mu)}{d\mu}d\mu$$

$$+\beta_R\int_{E_A} e^{i\varphi_R(\mu)x-\omega_0(\mu)t}G_{3,R}(\mu,t)\frac{d\varphi_R(\mu)}{d\mu}d\mu+\int_{-\infty}^{\infty}e^{i\lambda x-\omega_R(\lambda)t}\tilde{\hat{f}}_R(\lambda,\omega_R(\lambda),t)d\lambda\,,\ x>0\,,\ t>0\,.$$

Similarly,

$$(4.11)\quad 2\pi u^L(x,t)=\int_{-\infty}^{\infty}e^{i\lambda x-\omega_L(\lambda)t}\hat{u}_0^L(\lambda)d\lambda+i\int_{E_A}\varphi_L(-\mu)e^{i\varphi_L(-\mu)x-\omega_0(\mu)t}(\beta_L\varphi_L^2(-\mu)+\alpha_L)G_{0,L}(\mu,t)\left\{-\frac{d[\varphi_L(-\mu)]}{d\mu}\right\}d\mu$$

$$+\int_{E_A} e^{i\varphi_L(-\mu)x-\omega_0(\mu)t}(\beta_L\varphi_L^2(-\mu)+\alpha_L)G_{1,L}(\mu,t)\left\{-\frac{d[\varphi_L(-\mu)]}{d\mu}\right\}d\mu$$

$$-\beta_L i\int_{E_A}\varphi_L(-\mu)e^{i\varphi_L(-\mu)x-\omega_0(\mu)t}G_{2,L}(\mu,t)\left\{-\frac{d[\varphi_L(-\mu)]}{d\mu}\right\}d\mu$$

$$-\beta_L\int_{E_A} e^{i\varphi_L(-\mu)x-\omega_0(\mu)t}G_{3,L}(\mu,t)\left\{-\frac{d[\varphi_L(-\mu)]}{d\mu}\right\}d\mu+\int_{-\infty}^{\infty}e^{i\lambda x-\omega_L(\lambda)t}\tilde{\hat{f}}_L(\lambda,\omega_L(\lambda),t)d\lambda\,,\ x<0\,,\ t>0\,.$$

Thus, making $A$ larger, if necessary, and changing the variables of integration over $E_A$, by setting $\lambda=\varphi_R(\mu)$ and and $\lambda=-\varphi_L(-\mu)$, in (4.10) and (4.11), respectively, we obtain (4.5) and (4.6).

***Remark*** The representation formula (4.5) can also be written in the form

$$(4.12)\quad 2\pi u^R(x,t)=\int_{-\infty}^{\infty}e^{i\lambda x-\omega_R(\lambda)t}\hat{u}_0^R(\lambda)d\lambda-i\int_{\partial\mathcal{Z}_A}\lambda e^{i\lambda x-\omega_R(\lambda)t}(\beta_R\lambda^2+\alpha_R)G_{0,R}(\psi_R(\lambda),T)d\lambda$$

$$-\int_{\partial\mathcal{Z}_A}e^{i\lambda x-\omega_R(\lambda)t}(\beta_R\lambda^2+\alpha_R)G_{1,R}(\psi_R(\lambda),T)d\lambda+\beta_R i\int_{\partial\mathcal{Z}_A}\lambda e^{i\lambda x-\omega_R(\lambda)t}G_{2,R}(\psi_R(\lambda),T)d\lambda$$

$$+\beta_R\int_{\partial\mathcal{Z}_A}e^{i\lambda x-\omega_R(\lambda)t}G_{3,R}(\psi_R(\lambda),T)d\lambda+\int_{-\infty}^{\infty}e^{i\lambda x-\omega_R(\lambda)t}\tilde{\hat{f}}_R(\lambda,\omega_R(\lambda),t)d\lambda\,,$$

for fixed $T>0$, with $0<t\le T$ and $x>0$.

In particular, in the case $f_R\equiv 0$, the solution $u^R(x,t)$, for $(x,t)\in\mathbb{R}^+\times(0,T]$, is a combination of the exponentials $e^{i\lambda x-\omega_R(\lambda)t}$, $\lambda\in(\partial\mathcal{Z}_A^+)\cup\mathbb{R}$, expressed as an integral with appropriate measures, i.e., in Ehrenpreis form.

A formula analogous to (4.12) holds also for the solution $u^L(x,t)$.

***Note*** Analogous representations for the solution of interface problems for the pure biharmonic diffusion equation (this equation on the half-line has already been studied via the UTM in [5]) can be found as a straightforward corollary of the preceding exposition.